\documentclass{article}
\usepackage[margin=1in]{geometry}
\usepackage{amssymb}
\usepackage{amsmath}
\usepackage{amsthm}
\usepackage{hyperref}
\usepackage{etoolbox}
\usepackage{verbatim}

\newtheorem{theorem}{Theorem}[section]
\newtheorem*{theorem*}{Theorem}
\newtheorem{corollary}[theorem]{Corollary}

\newtheorem{lemma}[theorem]{Lemma}
\newtheorem{proposition}[theorem]{Proposition}

\theoremstyle{definition}{
\newtheorem{definition}[theorem]{Definition}
\newtheorem*{definition*}{Definition}

\newtheorem*{example*}{Example}
\newtheorem{remark}[theorem]{Remark}

}

\numberwithin{equation}{section}

\newcommand{\Imm}{\operatorname{Im}}
\newcommand{\gl}{\operatorname{GL}_2^+(\mathbb{Q})}

\newcommand{\slr}{\operatorname{SL}_2(\mathbb{R})}
\newcommand{\slz}{\operatorname{SL}_2(\mathbb{Z})}
\newcommand{\uh}{\mathcal{H}}
\newcommand{\suq}{\subseteq}
\begin{document}
\title{Ax-Lindemann and Andr\'e-Oort for a Nonholomorphic Modular Function}
\author{Haden Spence}





\maketitle
\begin{abstract}
The modular case of the Andr\'e-Oort Conjecture is a theorem of Andr\'e and Pila, having at its heart the well-known modular function $j$.  I give an overview of two other `nonclassical' classes of modular function, namely the quasimodular (QM) and almost holomorphic modular (AHM) functions.  These are perhaps less well-known than $j$, but have been studied by various authors including for example Masser, Shimura and Zagier.  It turns out to be sufficient to focus on a particular QM function $\chi$ and its dual AHM function $\chi^*$, since these (together with $j$) generate the relevant fields.  After discussing some of the properties of these functions, I go on to prove some Ax-Lindemann results about $\chi$ and $\chi^*$.  I then combine these with a fairly standard method of o-minimality and point counting to prove the central result of the paper; a natural analogue of the modular Andr\'e-Oort conjecture for the function $\chi^*$.
\end{abstract}
\bigskip
\section{Introduction}\label{sect:intro}
Let $\uh=\{\tau \in \mathbb{C}:\Imm\tau >0\}$ be the complex upper half plane.  We begin with the classical $j$-function, mapping from $\uh$ to $\mathbb{C}$, which is well-known to be a \emph{modular function}.  It is also transcendental, of course, but nonetheless has rich and beautiful arithmetic properties.  For any quadratic point $\tau\in\uh$, the number $j(\tau)$ is algebraic over $\mathbb{Q}$.  Such a point $j(\tau)$ is called a \emph{special point} or \emph{singular modulus}.  The singular moduli are precisely the $j$-invariants of elliptic curves with complex multiplication.  By a classical theorem of Schneider \cite{Schneider1957}, the imaginary quadratic numbers are the only algebraic elements of $\uh$ whose image under $j$ is also algebraic.

The special points of $j$ turn out to be a particular instance of a more general phenomenon.  A relation between coordinates in $\uh$ is called a \emph{geodesic relation} if it is of the form $\tau=c$ for some constant $c$, or $\tau_1=g\tau_2$ for some $g\in\gl$.  For each $N\in\mathbb{N}$, there is a polynomial $\Phi_N\in\mathbb{Z}[X,Y]$ with the property that
\[\Phi_N(j(\tau),j(g\tau))=0,\]
for all $\tau\in\uh$ and any $g\in\gl$ which, when written as a primitive integer matrix, has determinant $N$.  So we see that geodesic relations between coordinates $\tau_i\in\uh$ induce algebraic relations between their images $j(\tau_i)\in\mathbb{C}$.  In fact, geodesic relations are the only algebraic relations in $\uh$ that induce algebraic relations on the $j$-side.  This fact, known as the Ax-Lindemann theorem for $j$, was proven by Pila in \cite{Pila2011}; we will discuss it further later.  

We call subvarieties of $\mathbb{C}^n$ which arise in this way \emph{$j$-special subvarieties}.  So a $j$-special subvariety of $\mathbb{C}^n$ is (an irreducible component of) a variety cut out by some equations of the form $\Phi_N(z_i,z_j)=0$ and $z_i=j(\tau_i)$, for various values of $N$ and singular moduli $j(\tau_i)$.  In general, a \emph{$j$-special point} is a zero-dimensional $j$-special subvariety, that is, an $n$-tuple $(j(\tau_1),\dots,j(\tau_n))$ where every $\tau_i$ is a quadratic point.

The $j$-special subvarieties of $\mathbb{C}^n$ are rather sparse; given a random variety $V\suq\mathbb{C}^n$, we would not expect many $j$-special subvarieties to be contained within it.  Hence the following finiteness result, proven by Pila in \cite{Pila2011}.  We call the result Modular Andr\'e-Oort, since it is a special case of the full Andr\'e-Oort Conjecture, a statement about general Shimura varieties.  The full Andr\'e-Oort conjecture is known under GRH by work of Edixhoven, Klingler, Ullmo and Yafaev (see for instance \cite{Edixhoven2001}, \cite{Edixhoven2003}, \cite{Klingler2014} and \cite{Ullmo2006}), and is known unconditionally for $\mathcal{A}_g$, the moduli space of principally polarised abelian varieties of genus $g$; a result of Tsimerman, Pila et al: \cite{Pila2014}, \cite{Tsimerman2015}.  In turn, Andr\'e-Oort is a special case of the far-reaching Zilber-Pink conjecture, so Modular Andr\'e-Oort forms only a small part of a much larger picture.  There is a variety of literature on these topics; good starting points include surveys by Pila \cite{Pil} and Zannier \cite{Zannier2012}.  

\begin{theorem}[Andr\'e/Pila, Modular Andr\'e-Oort]\label{thrm:AOforj}
Let $V$ be a subvariety of $\mathbb{C}^n$.  Then $V$ contains only finitely many \emph{maximal} $j$-special subvarieties.
\end{theorem}
The ``maximal'' is certainly necessary; in general, a positive-dimensional $j$-special variety will always contain infinitely many proper $j$-special subvarieties.

\bigskip
The purpose of this paper is to prove an analogue of this result in a slightly different setting.  We will be discussing what happens when $j$ is replaced (or supplemented) by certain `modular-like' functions, which, while not modular functions in the strict sense, exhibit many of the same properties.  We will be focusing on two classes of near-modular function: the \emph{quasimodular functions} and the \emph{almost holomorphic modular functions}.  Specifically, we will look at a quasimodular function $\chi$ and a related almost holomorphic modular function $\chi^*$, defined as
\[\chi=1728\cdot\dfrac{E_2E_4E_6}{E_4^3-E_6^2},\qquad \chi^*=1728\cdot\dfrac{E_2^*E_4E_6}{E_4^3-E_6^2},\]
where $E_k$ is the usual $k$th Eisenstein series and
\[E_2^*(\tau)=E_2(\tau)-\dfrac{3}{\pi\Imm\tau}.\]
See Section \ref{sect:QMAHMFuncs} for details about the properties of $\chi$ and $\chi^*$.  One crucial fact is the existence of modular polynomials $\Psi_N\in\mathbb{Q}[X,Y,Z]$, having the property that
\[\Psi_N(\chi^*(g\tau),j(\tau),\chi^*(\tau))=0\]
for suitable $g\in\gl$ (compare with the classical modular polynomials $\Phi_N$).  In Section \ref{sect:Special} we construct these $\Psi_N$, using them and the $\Phi_N$ to construct what we call ``$(j,\chi^*)$-special'' varieties, directly analogous to the ``$j$-special'' varieties discussed above.  The $(j,\chi^*)$-special varieties exist only inside even Cartesian powers of $\mathbb{C}$; we consider $\mathbb{C}^{2n}$ as the Zariski closure of $\pi(\uh^n)$, where
\[\pi:\uh^n\to\mathbb{C}^{2n}\]
is defined by
\[\pi(\tau_1,\dots,\tau_n)= (j(\tau_1),\chi^*(\tau_1),\dots,j(\tau_n),\chi^*(\tau_n)).\]
The central theorem of this paper is the analogue of \ref{thrm:AOforj} in this setting:

\bigskip\noindent
\textbf{Theorem \ref{thrm:AOforChiStar2} (Andr\'e-Oort for $(j,\chi^*)$).}  \emph{Let $V\suq\mathbb{C}^{2n}$ be a variety.  Then $V$ contains only finitely many maximal $(j,\chi^*)$-special subvarieties.}

\bigskip
The proof is quite similar to that of \ref{thrm:AOforj}, following a standard strategy of o-minimality and point-counting developed by Pila and Zannier.  The majority of the novelty in its proof lies in the following ``Ax-Lindemann type'' result.  Loosely, it says that all the algebraic sets $S\suq \uh^n$ with $\pi(S)\suq V$ are accounted for by the ``weakly $\uh$-special varieties''.  These are defined in Section \ref{sect:Special}; put simply, they are subvarieties of $\uh^n$ cut out by geodesic relations.

\bigskip\noindent
\textbf{Corollary \ref{cor:ChiStarALReformed}.}  \emph{Let $V$ be an irreducible subvariety of $\mathbb{C}^{2n}$ and let $\mathcal{Z}=\pi^{-1}(V)\suq\uh^n$.  Then $\mathcal{Z}^\text{alg}$ is just the union of the weakly $\uh$-special subvarieties of $\mathcal{Z}$.}

\bigskip
Here, $\mathcal{Z}^\text{alg}$ is defined as the union of all connected, positive-dimensional, real semialgebraic subsets of $\mathcal{Z}$.  

It is important to emphasise the difficulties that lie in the proof of \ref{cor:ChiStarALReformed}.  Traditional Ax-Lindemann results have always relied heavily on the holomorphicity of the functions involved.  Since $\chi^*$ is not holomorphic, a lot of technical trickery is required to reach \ref{cor:ChiStarALReformed}.  By contrast, the argument needed to get from \ref{cor:ChiStarALReformed} to \ref{thrm:AOforChiStar2}, done in Section \ref{sect:AO}, is a standard ``point-counting'' argument applying the Pila-Wilkie theorem: a well known result from the theory of o-minimal structures.  

Along the way to our nonholomorphic Ax-Lindemann result, it turns out that we need an analogous Ax-Lindemann result for the quasimodular function $\chi$.  Such a result is of course interesting in its own right.  Since the quasimodular functions are holomorphic, not much is required for this.  We simply take an Ax-Lindemann result of Pila \cite{Pila2013}, concerning $j$ and its derivatives, and strengthen it slightly\footnote{In the presence of an Ax-Lindemann theorem for $\chi$, it is reasonable to ask whether there is also an Andr\'e-Oort result in that setting.  In fact it is not even obvious that we can formulate such a result.  The function $\chi$, unlike $\chi^*$, does not take algebraic values at quadratic points, so there is no clear notion of what a ``$(j,\chi)$-special'' variety should be.}.

The plan for this paper is as follows.  In Section 2, we discuss some of the basic background of quasimodular and almost holomorphic modular forms and functions.  In section 3, we discuss the special sets and special points of $\chi^*$ and $\chi$, which is of course crucial to any Andr\'e-Oort statement.  In section 4, the largest section of the paper, we prove the required Ax-Lindemann results, before concluding in section 5 with the proof of Theorem \ref{thrm:AOforChiStar2}.

\bigskip\noindent
\textbf{Note.} This work was undertaken during the course of the author's DPhil studies at the University of Oxford, and much of it is intended to appear in my DPhil thesis.

\bigskip\noindent
\textbf{Acknowledgements.} To my supervisor, Jonathan Pila, an enormous thank you; without Jonathan's unfailing support and guidance I would be utterly lost.  Jonathan has also provided many very helpful suggestions regarding the content and structure of this document specifically.  Thanks also go to Alan Lauder, whose suggestion for a miniproject got me interested in nonclassical modular functions in the first place, and to my father Derek, for proofreading and commenting on various early versions of this document.  I am pleased to thank the referee for their very thorough reading of this paper, including many helpful comments, corrections and suggestions.  This work was supported by the Engineering and Physical Sciences Research Council. 

\section{Quasimodular and Almost Holomorphic Modular Functions}\label{sect:QMAHMFuncs}
Let us begin by recalling some basic background about modular functions and Eisenstein series.
\begin{definition} 
 A modular function is a map $f:\uh\to\mathbb{C}$ with the following properties:
 \begin{itemize}  
  \item For any $\gamma\in\slz$ and any $\tau\in\uh$, we have $f(\gamma\cdot\tau)=f(\tau)$.  Here, as usual, elements of the group $\slr$ act on $\uh$ via M\"obius transformations,
  \[\begin{pmatrix}a&b\\c&d\end{pmatrix}\cdot\tau=\dfrac{a\tau+b}{c\tau+d}.\]
  \item $f$ is meromorphic on $\uh$.
  \item $f$ is ``meromorphic at $\infty$''.  That is, the Fourier expansion of $f$,
  \[f(\tau)=\sum_{k=-\infty}^{\infty}c_k\exp(2k\pi i\tau),\]
  has only finitely many negative terms.
 \end{itemize}
\end{definition}
One way to construct modular functions is through Eisenstein series.  The $k$th Eisenstein series $E_k$ is a function from $\uh$ to $\mathbb{C}$, defined as
\[E_k(\tau)=\dfrac{1}{2}\sum_{\substack{(m,n)\in\mathbb{Z}\\(m,n)=1}}\dfrac{1}{(m\tau+n)^k}.\]
For even $k\geq 4$, it is easy to see that $E_k$ converges absolutely, defining a holomorphic function, and further that
\[E_k(\gamma\tau)=(c\tau+d)^kE_k(\tau),\]
where $\gamma=\begin{pmatrix}a&b\\c&d\end{pmatrix}\in\slz$.  (For odd $k$, of course, the sum vanishes.)  So in particular the function
\[1728\cdot\dfrac{E_4^3}{E_4^3-E_6^2}\]
is invariant under the action of $\slz$; it turns out to be a modular function.  In fact this is simply the definition of the $j$-function.  It is well-known that the denominator 
\[\dfrac{1}{1728}(E_4^3-E_6^2),\]
which is known as the discriminant function and denoted $\Delta$, does not vanish anywhere on $\uh$, so $j$ is holomorphic on all of $\uh$.  It turns out that $j$ is really the only modular function we need to worry about, since the field of modular functions is just $\mathbb{C}(j)$.

\bigskip\noindent
\textbf{Note.} In the remainder of this section we will quite freely use facts proven in Zagier's excellent paper \cite[pages 18-22, 48-49, 58-60]{Zagier2008}.

\bigskip
So far we have only used the absolutely convergent Eisenstein series, namely those $E_k$ with $k\geq 4$.  The Eisenstein series $E_2$ does not converge absolutely, but by taking the terms of the sum in a suitable order, it does define a holomorphic function $E_2$.  It does not have the same transformation properties as the other $E_k$, but rather satisfies
\[E_2(\gamma\tau)=(c\tau+d)^2E_2(\tau)-\dfrac{6i}{\pi}c(c\tau+d).\]
Hence one can see that the modified function
\[E_2^*(\tau)=E_2(\tau)-\dfrac{3}{\pi\Imm\tau}\]
has the usual weight 2 transformation law, that is
\[E_2^*(\gamma\tau)=(c\tau+d)^2E_2^*(\tau).\]
The functions $E_2$ and $E_2^*$ are the prototype examples of, respectively, quasimodular forms and almost holomorphic modular forms.  
\begin{definition}\label{defn:AHMforms}
 A function $f:\uh\to\mathbb{C}$ is an \emph{almost holomorphic modular form} of weight $k$ if:
 \begin{itemize}
  \item $f(\tau)$ can be written as a polynomial in $(\Imm\tau)^{-1}$, with coefficients which are holomorphic functions, bounded as $\Imm\tau\to\infty$.
  \item $f$ satisfies the weight $k$ transformation law:
  \[f(\gamma\tau)=(c\tau+d)^kf(\tau).\]
 \end{itemize}
\end{definition}
\begin{definition}\label{defn:quasiforms}
 A function $f:\uh\to\mathbb{C}$ is a \emph{quasimodular form} of weight $k$ if it arises as the constant term (with respect to $(\Imm\tau)^{-1}$) of an almost holomorphic modular form of weight $k$.  Equivalently:
 \begin{itemize}
  \item $f(\tau)$ is a holomorphic function, bounded as $\Imm\tau\to\infty$.
  \item $f$ satisfies the modified transformation law:
  \[\dfrac{f(\gamma\tau)}{(c\tau+d)^k}=f(\tau)+\sum_{r=1}^p f_r(\tau)\left(\dfrac{c}{c\tau+d}\right)^r,\]
  for some holomorphic functions $f_r$, bounded as $\Imm\tau\to\infty$.
 \end{itemize}
\end{definition}
The graded algebra of almost holomorphic modular forms is generated over $\mathbb{C}$ by $E_2^*$, $E_4$ and $E_6$.  The graded algebra of quasimodular forms, similarly, is generated by $E_2$, $E_4$ and $E_6$.  In fact, these two graded algebras are isomorphic to each other via the map sending $E_2^*$ to $E_2$ and fixing $E_4$, $E_6$.  One can see this map as that sending an almost holomorphic modular form to its constant coefficient.  

For proofs of the various assertions made above, as well as more details about quasimodular and almost holomorphic modular forms in general, see \cite[pages 58-60]{Zagier2008}.  For this paper, we are more interested in quasimodular and almost holomorphic modular \emph{functions}.

\begin{definition}\label{defn:ahm+quasiFunctions}
 An almost holomorphic modular (or AHM) function is a quotient of almost holomorphic modular forms of the same weight.
 
 A quasimodular (or QM) function is a quotient of quasimodular forms of the same weight.
\end{definition}

The space of AHM functions and the space of QM functions are both obviously fields.  We will write $F^*$ for the field of AHM functions, and $\tilde{F}$ for the field of QM functions.  Each contains the field of classical modular functions.  These have been studied in a few places before, perhaps most notably by Masser in \cite[Appendix A]{Masser1975}.  Masser works with an AHM function he calls $\psi$, defined by $E_2^*E_4/E_6$.  This function has a singularity at $i$, so we work instead with a related function that has no singularities.

Define:
\[f=\dfrac{E_4E_6}{\Delta},\qquad\chi=E_2f,\qquad\chi^*=E_2^*f,\]
where $\Delta$ is again the discriminant function $(E_4^3-E_6^2)/1728$.  The function $f$ is then a meromorphic modular form of weight -2.  Since $\Delta$ does not vanish, none of these three functions have singularities inside $\uh$. 

Further, $\chi^*$ is an AHM function and $\chi$ is a QM function.  The function $\chi$ is holomorphic on $\uh$, but of course $\chi^*$ is only real analytic. 
We note for future use that
\[\chi^*(\tau)=\chi(\tau)-\dfrac{3}{\pi \Imm\tau}f(\tau),\]
and (by the transformation properties of $E_2$) that
\[\chi(\gamma\tau)=\chi(\tau)-\dfrac{6i}{\pi}\dfrac{c}{c\tau+d}f(\tau),\]
for all $\gamma=\begin{pmatrix}a&b\\c&d\end{pmatrix}\in\slz$.

\begin{proposition}\label{propn:fieldStructure}
 The fields $F^*$ and $\tilde{F}$ are characterised by:
 \[F^*=\mathbb{C}(j,\chi^*),\qquad\tilde{F}=\mathbb{C}(j,\chi).\]
 Moreover, $F^*$ and $\tilde{F}$ are isomorphic via the map fixing $j$ and sending $\chi^*$ to $\chi$.
 \end{proposition}
 \begin{proof}[Proof Sketch]
  Zagier proves in \cite[Proposition 20, page 59]{Zagier2008} that the graded algebras of QM and AHM forms are generated by $E_4$, $E_6$ and (respectively) $E_2$ or $E_2^*$.  Given that fact, it is a simple exercise to write down a generating set for the `monomial quotients' of QM and AHM forms, and see that they are all expressible as rational functions of $j$ and $\chi$ or $\chi^*$.
  
  The isomorphism of fields is induced directly by the isomorphism between the graded algebras of QM and AHM forms.
 \end{proof}

The following will also be of use.
\begin{theorem}\label{thrm:algIndependence}
 The functions $j$, $\chi$ and $f$ are algebraically independent over $\mathbb{C}$.
 \end{theorem}
 \begin{proof}
  Follows easily from the standard fact that $j$, $j'$ and $j''$ are algebraically independent functions over $\mathbb{C}$.  See for instance Zagier \cite[page 49]{Zagier2008}.
 \end{proof}
Our intent is to discuss the special sets corresponding to the functions $\chi^*$ and $\chi$.  Such things do exist; they are the subject of the next section.
\section{Special Sets}\label{sect:Special}
\subsection{New Modular Polynomials}
Our discussion of special sets begins with the following proposition involving the construction of some modular polynomials for $\chi^*$.  Although this follows fairly easily from facts known about $j$ and its derivatives, together with the upcoming Lemma \ref{lma:isomorphism}, the explicit existence of these polynomials seems not to have been noted before. The construction is very similar to the standard construction of the usual modular polynomials; we follow Zagier \cite[Proposition 23, pages 68-69]{Zagier2008} closely.
\begin{proposition}\label{propn:ModRelations}
 For a positive integer $N$, let $M_N$ be the set of primitive integer matrices $g\in\gl$ with determinant $N$.  For each such $N$, there is a nonzero polynomial $\Psi_N\in\mathbb{Q}[X,Y,Z]$, irreducible over $\mathbb{C}$, such that 
 \[\Psi_N(\chi^*(g\tau),j(\tau),\chi^*(\tau))=0\]
 for each $g\in M_N$ and all $\tau\in\uh$.
\end{proposition}
 \begin{proof}
  The set \[D_N=\left\{\begin{pmatrix}a&b\\0&d\end{pmatrix}:a,b,d\in\mathbb{N}, ad=N, 0\leq b<d, \gcd(a,b,d)=1\right\}\]
  is a full set of representatives for $M_N$ under the action of $\slz$.  That is, for all $g\in M_N$ there is some $g'\in D_N$ and $\gamma\in\slz$ such that $\gamma g'=g$.  (This is a standard fact; see for instance Lang \cite{Lang1987} or Diamond/Shurman \cite[Exercise 1.2.11]{Diamond2005}.)
  
  We will consider a polynomial in $X$, defined by
    \begin{equation}\label{ComboPoly}
      \prod_{g\in D_N}(X-\chi^*(g\tau)).
    \end{equation}
  Clearly (for each $\tau$) this is 0 if and only if $X$ is $\chi^*(h\tau)$, for some $h\in D_N$.  Thanks to the invariance of $\chi^*$ under $\slz$, this holds if and only if $X$ is $\chi^*(h\tau)$ for some $h\in M_N$.

  Let $\gamma\in\slz$.  For each $g\in D_N$, we have $g\cdot\gamma=\gamma'\cdot h$, for some other $\gamma'\in\slz$ and some $h\in D_N$.  So by the invariance of $\chi^*$, we have
  \[\chi^*(g\cdot\gamma\tau)=\chi^*(\gamma'\cdot h\tau)=\chi^*(h\tau).\]
  Thus the map $\tau\mapsto\gamma\tau$ induces a permutation of the set
  \[S_N=\{\chi^*(g\tau):g\in D_N\}.\]
  In fact, the described action of $\slz$ on $S_N$ is transitive.  Indeed, any $g\in D_N$ can be written as
  \[g=\gamma h\gamma',\qquad \gamma,\gamma'\in\operatorname{GL}_2(\mathbb{Z}),\]
  with $h$ in Smith Normal Form, meaning it is a diagonal matrix $\begin{pmatrix}A&0\\0&D\end{pmatrix}$, with $A|D$ (see for instance \cite[Exercise 19, page 470]{Dummit2004}).  Further, $h$ must be primitive since $g$ is, whence $A=1$ and $D=N$.  By replacing $\gamma,\gamma'$ by $\gamma\begin{pmatrix}1&0\\0&-1\end{pmatrix}$ and $\begin{pmatrix}1&0\\0&-1\end{pmatrix}\gamma'$ if necessary, we can ensure that they are in fact elements of $\slz$.  The claimed transitivity follows immediately\footnote{I thank David Speyer for showing me the proof of this fact, which is taken as read in many texts.}.
   
  Each coefficient of $X$ in the polynomial (\ref{ComboPoly}) is a symmetric polynomial in the functions $\chi^*(g\tau)$, $g\in D_N$, so each coefficient must be invariant under $\tau\mapsto\gamma\tau$.  Moreover, if $g=\begin{pmatrix}a&b\\0&d\end{pmatrix}\in D_N$, then
  \[\Imm(g\tau)=\dfrac{a\Imm\tau}{d}.\]
  Hence each coefficient is a polynomial in $1/\Imm\tau$ with coefficients which are meromorphic functions on $\uh$.  Since they are also $\slz$-invariant, each coefficient is therefore an element of the field of AHM functions $F^*=\mathbb{C}(j,\chi^*)$, so can be written as a quotient of complex polynomials in $j$ and $\chi^*$.  In each such rational function, we can replace instances of $j$ and $\chi^*$ with variables $Y$ and $Z$.  If we do this for each coefficient, we get a polynomial
  \[\Psi_N^0(X,Y,Z)\in\mathbb{C}(Y,Z)[X]\]
  with $\Psi_N^0(X,j(\tau),\chi^*(\tau))=0$ if and only if $X=\chi^*(g\tau)$ for some $g\in M_N$.
  
  Next, note that
  \[\chi^*(x+iy)=\dfrac{E_2E_4E_6}{\Delta}-\dfrac{3}{\pi y}\cdot\dfrac{E_4E_6}{\Delta}.\]
  Each of the Eisenstein series and $\Delta$ has a power series expansion in $q=e^{2\pi i z}$, with integer coefficients.  The coefficient of the leading term in each case is 1; the coefficients of the $q$-expansions of $E_2$, $E_4$ and $E_6$ are given, for example, in \cite[pages 17 and 19]{Zagier2008}, and the $q$-expansion of $\Delta$ is easily calculated from those.
  
  Hence $\chi^*(x+iy)$ is a polynomial in $3/\pi y$ with coefficients that are Laurent series in $q$ with integer coefficients and leading term $q^{-1}$.  The function $j$ also has a $q$-expansion, which is just an integer Laurent series in $q$, again with leading term $q^{-1}$.  We will use this to show that $\Psi_N^0$ is defined over $\mathbb{Q}$. 
  
  We have (writing $y=\Imm\tau$)
  \begin{align*}
   \Psi_N^0(X,j(\tau),\chi^*(\tau))&=\prod_{\substack{ad=N\\d>0}}\prod_{\substack {0\leq b<d\\(a,b,d)=1}}\left(X-\chi^*\left(\dfrac{a\tau+b}{d}\right)\right)\\
   &=\prod_{\substack{ad=N\\d>0}}\prod_{\substack {0\leq b<d\\(a,b,d)=1}}\left(X-\sum_{n=-1}^{\infty} c_n \zeta_d^{nb}q^{na/d}+\dfrac{d}{a}\dfrac{3}{\pi y}\sum_{n=-1}^{\infty} c_n' \zeta_d^{nb}q^{na/d}\right),
  \end{align*}
  where $\zeta_d=e^{2\pi i/d}$ and $c_n,c_n'\in\mathbb{Z}$.  The inner product is a polynomial in $3/\pi y$, with coefficients which are Laurent series in $q^{1/d}$ with coefficients from $\mathbb{Z}\left[\frac{d}{a},\zeta_d\right]$, and leading term no smaller than $-a/d$.  But it is 1-periodic, so the fractional powers of $q$ must cancel out.  Further, the coefficients in the resulting $q$-expansions must be in $\mathbb{Z}[\frac{d}{a}]$, since every Galois conjugation $\zeta_d\mapsto \zeta_d^r$, where $r\in(\mathbb{Z}/d\mathbb{Z})^*$, fixes the inner product; the numbers $b$ and $rb$ range over the same set.
  
  So each coefficient $f_k$ of $X^k$ in $\Psi_N^0$ is a polynomial in $3/\pi y$ with coefficients which are rational Laurent series in $q$.  Each coefficient is also equal to a quotient of polynomials $p_k$ and $q_k$ in $j$ and $\chi^*$, thus
  \[p_k(j,\chi^*)=f_k\cdot q_k(j,\chi^*).\]
  If we compare the coefficients of $(3/\pi y)^k$ on each side, we get various equalities between $q$-expansions.  The coefficients of those $q$-expansions are $\mathbb{Q}$-linear in the coefficients of $p_k$ and $q_k$.  So we get a homogeneous system of $\mathbb{Q}$-linear equations holding for the coefficients of $p_k$ and $q_k$.  This system certainly has a solution since $p_k$ and $q_k$ exist.  By basic linear algebra, the solution can be chosen to be rational up to scaling, ie. $p_k$ and $q_k$ are in $\lambda\mathbb{Q}[Y,Z]$, for some $\lambda$.  In particular, $p_k/q_k$ can be rewritten as a quotient of rational polynomials.
  
  Thus $\Psi_N^0\in\mathbb{Q}(Y,Z)[X]$.  Finally, since, as noted earlier, $\slz$ acts transitively on $S_N$, no subproduct of
  \begin{equation*}\label{eqn:ModPolyProduct}\prod_{g\in D_N}(X-\chi^*\circ g)=\Psi_N^0(X,j,\chi^*)\in F^*[X]\end{equation*}
  can have coefficients that are $\slz$-invariant.  Hence $\Psi_N^0(X,j,\chi^*)$ is irreducible over $F^*$.  In particular, $\Psi_N^0(X,Y,Z)$ is irreducible over $\mathbb{C}(Y,Z)$ as a polynomial in $X$.  It is also monic in $X$, so if we clear the denominators in $Y$ and $Z$ exactly, we get an irreducible polynomial $\Psi_N\in\mathbb{Q}[X,Y,Z]$ having the required properties.
 \end{proof}
In the above, we have made essential use of the fact that $M_N$ is represented (up to the action of $\slz$) by the finitely many \emph{upper triangular} matrices in $D_N$.  Since $\chi^*$ is $\slz$-invariant, it is enough that the relation
\[\Psi_N(\chi^*(g\tau),j(\tau),\chi^*(\tau))=0\]
holds for $g\in D_N$; that implies the relation for all of $M_N$.  This is not the case for the QM function $\chi$, which only exhibits nice properties with respect to upper triangular elements of $\gl$.  The best we can do is the following.
\begin{lemma}\label{lma:isomorphism}
 Let $g_1\dots,g_k$ be \emph{upper triangular} primitive integer matrices with positive determinant, and consider the two fields:
 \[A=F^*(j\circ g_1,\dots,j\circ g_k,\chi^*\circ g_1,\dots,\chi^*\circ g_k)\]
 and
 \[B=\tilde{F}(j\circ g_1,\dots,j\circ g_k,\chi\circ g_1,\dots,\chi\circ g_k),\]
 considered as fields of real analytic functions, defined locally.
 Then $A$ and $B$ are isomorphic via the map
 \[\chi^*\mapsto \chi,\qquad \chi^*\circ g_i\mapsto\chi\circ g_i,\]
 fixing $j$ and all of the $j\circ g_i$.
 \end{lemma}
 \begin{proof}   
  The map is clearly a well-defined bijection.  If some $\chi^*\circ g_i$ and $j\circ g_i$ satisfy a polynomial equation $p(\chi^*\circ g_1,j\circ g_1,\dots,\chi^*\circ g_k,j\circ g_k)=0$, then (by comparing growth rates) every coefficient of $1/\Imm\tau$ on the left hand side must vanish.  In particular, the constant term $p(\chi\circ g_1,j\circ g_1,\dots,\chi\circ g_k,j\circ g_k)$ must vanish.  That is, the same polynomial equation holds for the $\chi\circ g_i$ and $j\circ g_i$, so the map is indeed an isomorphism.
 \end{proof}
It follows that
\[\Psi_N(\chi(g\tau),j(\tau),\chi(\tau))=0,\] for all \emph{upper-triangular} primitive integer matrices of determinant $N$.  The relation fails for any matrix which is not upper triangular; simply look at the transformation law satisfied by $\chi$.  

The existence of the polynomials $\Psi_N$ allows us to say the same thing about $\chi^*$ that is true of $j$, namely: \emph{geodesic relations between coordinates $\tau_i\in\uh$ induce algebraic relations between their images $j(\tau_i),\chi^*(\tau_i)\in\mathbb{C}$}.  Similarly, we can say of $\chi$ that \emph{upper triangular geodesic relations} 
induce algebraic relations.  Hence, just as for $j$, we can talk about the special varieties of $\mathbb{C}^n$ corresponding to $\chi^*$ and $\chi$. 
\subsection{Special Subvarieties}
We will discuss various related types of special subvariety:
\begin{itemize}
 \item $\uh$-special and weakly $\uh$-special subvarieties of $\uh^n$.
 \item $j$-special and weakly $j$-special subvarieties of $\mathbb{C}^n$.
 \item $(j,\chi^*)$-special, weakly $(j,\chi^*)$-special and weakly $(j,\chi)$-special subvarieties of $\mathbb{C}^{2n}$.
 \item $\chi^*$-special, weakly $\chi^*$-special and weakly $\chi$-special subvarieties of $\mathbb{C}^n$.
\end{itemize}
We shall see that the weakly $\chi^*$-special and weakly $\chi$-special subvarieties turn out to be very similar objects.  However, the naive approach towards producing a ``truly $\chi$-special'' subvariety does not appear to work.  The same goes for $(j,\chi)$-special subvarieties.

\subsubsection{\texorpdfstring{$\uh$}{H}-special and \texorpdfstring{$j$}{j}-special Varieties}
We will start with the $\uh$-special subvarieties.
\begin{definition}
Let $n\in\mathbb{N}$.

 Let $S_0\cup S_1\cup\dots\cup S_k$ be a partition of $\{1,\dots,n\}$, where $k\geq 0$ and $S_i\ne\emptyset$ for $i>0$.
 For each $s\in S_0$, choose any point $q_s\in\uh$.  For each $i>0$, let $s_i$ be the least element of $S_i$ and for each $s_i\ne s\in S_i$ choose a geodesic matrix $g_{i,s}\in\gl$.  A \emph{weakly $\uh$-special} subvariety of $\uh^n$ is a set of the form
 \begin{equation*}
   \bigl\{(\tau_1,\dots,\tau_n)\in\uh^n:\tau_s=q_s\text{ for }s\in S_0, \tau_{s}=g_{i,s}\tau_{s_i} \text{ for }s\in S_i, s\ne s_i,i=1,\dots,k\bigr\},
 \end{equation*}
 for some given data $S_i$, $q_s$, $g_{i,s}$.
 
 A weakly $\uh$-special subvariety is \emph{$\uh$-special} if the constant factors $q_s$ are imaginary quadratic numbers for all $s\in S_0$. 
\end{definition}
\begin{remark}
 What we call a ``weakly $\uh$-special variety'' is elsewhere in the literature referred to as a ``geodesic variety''.  We have chosen our terminology differently here so that it meshes as closely as possible with the terminology we use for $j$-special varieties and so forth.
\end{remark}
This definition may look more complex than it actually is.  Put more loosely, a weakly $\uh$-special variety is simply one defined by some number of equations of the form $\tau_i=q_i$ or $\tau_i=g_{ij}\tau_j$, with $q_i$ constants and $g_{ij}\in\gl$.  If it happens that the $q_i$ are quadratic points then the variety is $\uh$-special.  

This theme continues for all the types of special variety we will define here; a special variety will be defined by some combination of:
\begin{itemize}
 \item Equations requiring some coordinate to be constant (perhaps a ``special'' constant).
 \item Equations coming from the modular polynomials.
\end{itemize}
These other types of special variety will all arise as (the Zariski closures of) the images of $\uh$-special varieties under various maps.  The easiest to deal with are the well-known $j$-special varieties.
\begin{definition}\label{defn:specSub}
 Let $n\in\mathbb{N}$ and let $S_0\cup S_1\cup\dots\cup S_k$ be a partition of $\{1,\dots,n\}$, where $k\geq 0$ and $S_i\ne\emptyset$ for $i>0$.
  For each $s\in S_0$, choose a point $j_s\in\mathbb{C}$.  For each $i>0$, let $s_i$ be the least element of $S_i$ and for each $s_i\ne s\in S_i$ choose a positive integer $N_{i,s}$.  A weakly $j$-special subvariety of $\mathbb{C}^{n}$ is an irreducible component of a subvariety of the form
  \begin{equation*}\{(z_1,\dots,z_n)\in\mathbb{C}^n:z_s=j_s\text{ for }s\in S_0,\Phi_{N_{i,s}}(z_{s_i},z_s)=0\text{ for }s\in S_i, s\ne s_i,i=1,\dots,k\}\end{equation*}
  for some given data $S_i$, $j_s$, $N_{i,s}$.
  
  A weakly $j$-special variety is $j$-special if all of the constant factors $j_s$ are singular moduli, ie. $j_s=j(\tau_s)$ for some quadratic $\tau_s\in\uh$.
\end{definition}
The $j$-special (resp. weakly $j$-special) varieties of $\mathbb{C}^n$ are precisely those varieties that arise as the image, under the map 
\[(\tau_1,\dots,\tau_n)\mapsto(j(\tau_1),\dots,j(\tau_n)),\]
of an $\uh$-special (resp. weakly $\uh$-special) subvariety of $\uh^n$.  These are the special varieties discussed in \ref{thrm:AOforj}.

\subsubsection{\texorpdfstring{$(j,\chi^*)$}{(j,Chi-Star)}-special Varieties}\label{sect:jChiSpecials}
The $(j,\chi^*)$-special subvarieties are slightly more intricate.  We start in the simplest positive dimensional case, considering the variety
\begin{equation*}V_N'=\big\{(W,X,Y,Z)\in\mathbb{C}^4:\Phi_N(W,Y)=0,\,\Psi_N(X,Y,Z)=0,\Psi_N(Z,W,X)=0\big\}\suq\mathbb{C}^4.\end{equation*}
By counting equations, $\dim_{\mathbb{C}}V_N'$ is at most 2.  In fact, $\dim_{\mathbb{C}}V_N'=2$.  To see this, note that $V_N'$ contains the set
\[S_g=\{(j(\tau),\chi^*(\tau),j(g\tau),\chi^*(g\tau)):\tau\in\uh\}\]
for any $g\in M_N$.  Since $j$ and $\chi^*$ are algebraically independent, $S_g$ cannot be contained in any algebraic curve; hence $\dim_\mathbb{C}V_N'>1$.

I believe that the variety $V_N'$ is always irreducible.  I have calculated the first few modular polynomials $\Psi_N$ to see in each case that $V_N'$ is irreducible, but so far have been unable to prove that this is the case for all $N$.  I leave this as an open problem, which fortunately has no impact whatsoever on the wider picture: by real analytic continuation, $V_N'$ has an irreducible component containing $S_g$.  Call this component $V_N$; it is still 2-dimensional.  Since it contains $S_g$, the variety $V_N$ in fact contains all the $S_g$, $g\in M_N$, by modularity of $j$ and $\chi^*$.  Moreover, by \ref{lma:isomorphism}, $V_N$ contains
\[S_g'=\{(j(\tau),\chi(\tau),j(g\tau),\chi(g\tau)):\tau\in\uh\}\]
for any upper triangular $g\in M_N$.  These $V_N$ will form the building blocks of $(j,\chi^*)$-special varieties.
\begin{definition}
  Let $n\in\mathbb{N}$ and let $S_0\cup S_1\cup\dots\cup S_k$ be a partition of $\{1,\dots,n\}$, where $k\geq 0$ and $S_i\ne\emptyset$ for $i>0$.
  For each $s\in S_0$, choose $\tau_s\in\uh$ and let $(j_s,c_s)=(j(\tau_s),\chi^*(\tau_s))\in\mathbb{C}^2$.  For each $i>0$, let $s_i$ be the least element of $S_i$ and for each $s_i\ne s\in S_i$ choose a positive integer $N_{i,s}$.  A weakly $(j,\chi^*)$-special subvariety of $\mathbb{C}^{2n}$ is an irreducible component of a subvariety of the form
  \begin{multline*}
   \bigl\{(w_1,z_1,\dots,w_n,z_n)\in\mathbb{C}^{2n}:(w_s,z_s)=(j_s,c_s)\text{ for }s\in S_0,\\ (w_s,z_s,w_{s_i},z_{s_i})\in V_{N_{i,s}} \text{ for }s\in S_i, s\ne s_i,i=1,\dots,k\bigr\},
  \end{multline*}
  for some given data $S_i$, $(j_s,c_s)$, $N_{i,s}$.
  
  A weakly $(j,\chi^*)$-special variety is $(j,\chi^*)$-special if every constant factor $(j_s,c_s)$ is of the form $(j(\tau_s),\chi^*(\tau_s))$ for some quadratic point $\tau_s\in\uh$.
\end{definition}
Every weakly $(j,\chi^*)$-special variety arises as the Zariski closure of the image of a weakly $\uh$-special variety under the map 
\[(\tau_1,\dots,\tau_n)\mapsto(j(\tau_1),\chi^*(\tau_1),\dots,j(\tau_n),\chi^*(\tau_n)).\]
One can see from the polynomials $\Psi_N$ that $\chi^*(\tau)$ is an algebraic number whenever $\tau\in\uh$ is quadratic; this also follows from the work of Masser \cite[Appendix A]{Masser1975}.  This is why $(j,\chi^*)$-special varieties are called such; all their constant factors are special algebraic numbers. 

The weakly $(j,\chi)$-special varieties differ from $(j,\chi^*)$-special varieties only in their constant factors.  The definition is identical, except that the constant factors $(j_s,c_s)$ are chosen to be of the form $(j(\tau),\chi(\tau))$.   Unlike in the AHM case, these $(j,\chi)$-special varieties do not arise as the Zariski closure of any arbitrary $\uh$-special set.  
\begin{definition}\label{defn:GUTvar}
A geodesic upper triangular (or GUT) variety is a weakly $\uh$-special variety for which all of the relations $g_{i,s}$ are upper triangular matrices.
\end{definition}
Since $\chi$ behaves nicely only under upper-triangular matrices, one can see that weakly $(j,\chi)$-special varieties arise only from GUT varieties.  A weakly $(j,\chi)$-special variety is the Zariski closure of the image of a GUT variety under the map
\[(\tau_1,\dots,\tau_n)\mapsto(j(\tau_1),\chi(\tau_1),\dots,j(\tau_n),\chi(\tau_n)).\]
\begin{remark}
 The polynomials $\Psi_N$ only work for $\chi$ and $\chi\circ g$ when all of the relevant matrices are upper triangular.  As a result, numbers $\chi(\tau)$ with $\tau$ quadratic are \emph{not} algebraic in general.  Diaz has proven and conjectured various results about these points and their transcendence properties in his paper \cite{Diaz2000}.  So $\chi$ seems not to have any points that we could reasonably call special points, other than perhaps the $\slz$-translates of $i$ and $e^{2\pi i/3}$, where $\chi$ vanishes. 
 
 This is why we have not attempted to define a notion of $(j,\chi)$-special variety; the naive approach does not seem to produce a correct definition and it is not immediately clear whether a correct such notion exists.  This is also why, in \ref{defn:GUTvar}, we have not defined any notion of ``$\uh$-special GUT variety''.  We might, for instance, have required all of the constant factors to be quadratic, or to be special in some other way, but this seems unlikely to produce a meaningful outcome since $\chi$ apparently has no special points.  
\end{remark}

\subsubsection{\texorpdfstring{$\chi^*$}{Chi-Star}-special Varieties}
The final special varieties we are interested in are the $\chi^*$-special varieties.  The idea is that, while any individual polynomial $\Psi_N$ introduces a dependence on a $j$-coordinate, multiple different relations induced by $\Psi_N$ can cancel each other out and introduce a relation that exists between the $\chi^*\circ g_i$ alone, not involving $j$.

Unfortunately, it seems difficult to isolate the specific polynomials that arise in this process.  So the easiest definition to use for $\chi^*$-special varieties is the following:
\begin{definition}
 A (weakly) $\chi^*$-special variety is an irreducible component of the Zariski closure of the projection of a (weakly) $(j,\chi^*)$-special variety onto the coordinates corresponding to $\chi^*$-variables.  
\end{definition}
Equivalently, a $\chi^*$-special (resp. weakly $\chi^*$-special) variety is an irreducible component of the Zariski closure of the image of an $\uh$-special (resp. weakly $\uh$-special) set under the map
\[(\tau_1,\dots,\tau_n)\mapsto(\chi^*(\tau_1),\dots,\chi^*(\tau_n)).\]
As before, the weakly $\chi$-special varieties differ from the weakly $\chi^*$-special varieties only in their constant factors, and there is no obvious concept of a $\chi$-special variety.

In what follows, we prove a few results about the possible shape of weakly $\chi^*$-special varieties.

\begin{proposition}\label{propn:infiniteGs}
 Let $N\geq 2$ and let $S=\chi^{-1}\{0\}$.  Then there is an upper triangular $g\in M_N$ such that the set
 \[\{\chi(gs):s\in S\}\]
 is infinite.
 \end{proposition}
 \begin{proof}
  For any $\tau\in\uh$ which is $\slz$-equivalent to $i$, the Eisenstein series $E_6$ is equal to 0.  In particular, $\slz\cdot i\suq S$.  So we only need to show that (for some $g$) $\chi(g(\gamma \cdot i))$ takes infinitely many values as $\gamma$ varies.
 This is easy to see simply by considering matrices of the form
 \[g=\begin{pmatrix}N&0\\0&1\end{pmatrix}\in D_N,\qquad \gamma_n=\begin{pmatrix}1&-1\\1-nN&nN\end{pmatrix}\in\slz.\]
 Then we get
 \[g\cdot\gamma_n=\begin{pmatrix}N&-1\\1-nN&n\end{pmatrix}\cdot\begin{pmatrix}1&0\\0&N\end{pmatrix},\]
 so using the transformation law for $\chi$, we have
 \[\chi(g(\gamma_n\cdot \tau))=\chi(\tau/N)-\dfrac{6i}{\pi}\dfrac{1-nN}{(1-nN)(\tau/N)+n}f(\tau/N),\]
 where $f=E_4E_6/\Delta$.  Setting $\tau=i$, the above expression clearly takes infinitely many values as $n$ varies, provided that $f(i/N)\ne0$, which is clear since the only zeros of $E_4$ and $E_6$ are $\slz$-equivalent to $i$ or $\rho$. 
 \end{proof}
\begin{corollary}\label{cor:jPolys}
The modular polynomial $\Psi_N(X,Y,Z)$ is nonconstant in $Y$ for all $N\geq 2$.
\end{corollary}
\begin{proof}
 Otherwise $\Psi_N=\Psi_N(X,Z)$, and then by \ref{propn:infiniteGs}, $\Psi_N(0,Z)$ has infinitely many solutions, and so is identically zero.  Since $\Psi_N$ is irreducible, this would mean that $\Psi_N(X,Z)$ is a constant multiple of $X$, which is clearly false.
\end{proof}
\begin{remark}
 The modular polynomial $\Psi_1(X,Y,Z)$ is just $X-Z$.  By the above it is the only modular polynomial which does not depend on $Y$.
\end{remark}

\begin{corollary}
 Let $n\geq 1$.  For each $1\leq i\leq n$, let $g_i$ be an \emph{upper triangular} primitive integer matrix with positive determinant $N_i$.  Suppose that not all the $N_i$ are equal to 1.  Then the Zariski closure of the set
 \[S=\{(\chi(\tau),\chi(g_1\tau),\dots,\chi(g_n\tau)):\tau\in\uh\}\]
 has complex dimension exactly 2.
 Similarly, the Zariski closure of
 \[S^*=\{(\chi^*(\tau),\chi^*(g_1\tau),\dots,\chi^*(g_n\tau)):\tau\in\uh\}\]
 has complex dimension exactly 2.
 \end{corollary}
 \begin{proof}
  Consider the $(j,\chi^*)$-special variety $W\suq\mathbb{C}^{2n+2}$, defined by
  \[W=\{(X_0,\dots,X_n,Y_0,\dots,Y_n):(X_0,Y_0,X_i,Y_i)\in V_{N_i}, 1\leq i \leq n\}.\]
  $W$ is a 2-dimensional variety and contains the sets
  \[\{(j(\tau),j(g_1\tau),\dots,j(g_n\tau),\chi^*(\tau),\chi^*(g_1\tau),\dots,\chi^*(g_n\tau))\} \]
  and
  \[\{(j(\tau),j(g_1\tau),\dots,j(g_n\tau),\chi(\tau),\chi(g_1\tau),\dots,\chi(g_n\tau))\}.\]
  So the sets $S$ and $S^*$ are each contained in the projection of $W$ onto the $Y_i$-coordinates (which correspond to $\chi^*$-variables). Since $\dim W=2$, the projection certainly has dimension at most 2.  So we need to show that $S$ is not contained in a curve, which is immediate from \ref{propn:infiniteGs} since not all the $N_i$ are equal to 1.  By \ref{lma:isomorphism}, $S^*$ cannot be contained in any curve either.
 \end{proof}
\begin{corollary}\label{cor:ChiSpecialSubsOfC2}
 The only positive-dimensional, proper weakly $\chi^*$-special (or indeed weakly $\chi$-special) subvarieties of $\mathbb{C}^2$ are the diagonal $X=Y$ and the horizontal and vertical lines.
 \end{corollary}
 \begin{proof}
  Immediate.
 \end{proof}
\section{Ax-Lindemann Theorems}
\subsection{The Pila-Wilkie Counting Theorem}
In the remainder of this document we will use, several times, the theory of o-minimal structures, a branch of model theory.  The study of o-minimal structures has been used to great success by Habegger, Masser, Pila, Tsimerman, Zannier and many others to work with problems in functional transcendence, diophantine geometry and other areas.  Readers unfamiliar with the topic can see the excellent book by van den Dries \cite{Dries1998} or surveys by Pila \cite{Pil}.  

The crucial theorem we need for the purposes of this article is the following.
\begin{theorem}[Pila-Wilkie Counting Theorem]\label{thrm:PilaWilkie}
 Let $Z\suq \mathbb{R}^n$ be a definable set in an o-minimal structure $(\mathbb{R},\{+,\cdot,\dots\},\{<,\dots\})$.
 
 For each $\epsilon>0$ and each $k\in\mathbb{N}$ there is a constant $c(Z,k,\epsilon)$, such that for every $T\in\mathbb{N}$, we have
 \begin{equation*}\#\Bigl\{(x_1,\dots,x_n)\in \overline{\mathbb{Q}}^n\cap Z\setminus Z^\text{alg}:\max_i[\mathbb{Q}(x_i):\mathbb{Q}]\leq k, \max_i \operatorname{Ht}(x_i)\leq T\Bigr\}\leq c(Z,k,\epsilon)T^\epsilon.\end{equation*}
\end{theorem}
This may require some explanation.  The set $Z$ here is supposed to be ``definable in an o-minimal structure''.  For details about what this means, one can see the surveys cited above.  It is enough to know that there is a certain class of subsets of $\mathbb{C}^n=\mathbb{R}^{2n}$ which will be called ``definable in the o-minimal structure $\mathbb{R}_\text{an,exp}$,'' or just ``definable''.  Crucially, the graphs of the functions $j$, $\chi$ and $\chi^*$, restricted to any $\slz$-translate of the standard fundamental domain 
\[\mathbb{D}=\left\{\tau\in\uh:-\frac{1}{2}<\operatorname{Re}\tau\leq\frac{1}{2},|\tau|>1\right\},\] 
are all definable sets.  This follows from the fact that each of the functions can be written as a sum of convergent $q$-expansions, but can also be seen using a result of Peterzil and Starchenko \cite{Peterzil2004} about the definability of the Weierstrass $\wp$-function, going via the theory of elliptic curves.

Consequently, for any variety $V\suq\mathbb{C}^{kn}$, the restricted preimage
\[\pi^{-1}(V)\cap\mathbb{D}^n\]
is a definable set whenever $\pi:\uh^n\to\mathbb{C}^{kn}$ is some combination of the maps $j$, $\chi$ and $\chi^*$.  We will be applying \ref{thrm:PilaWilkie} to sets of precisely this form.

Given a definable set $Z$, we can consider what is called the ``transcendental part of $Z$'', $Z\setminus Z^\text{alg}$, where $Z^\text{alg}$ is the union of all connected, positive-dimensional, real semialgebraic subsets of $Z$.  Pila-Wilkie tells us that the number of algebraic points in $Z\setminus Z^\text{alg}$, of degree less than some fixed $k$ and height at most $T$, grows more slowly than any positive power of $T$.
Hence, if we can prove that a given definable set $Z$ contains `too many' algebraic points of a given height and degree, then $Z$ must contain an arc of a real algebraic curve.  

If $Z=\pi^{-1}(V)\cap\mathbb{D}^n$ as above, our next task is to find out which real algebraic curves can exist within such preimages.  For this we need what is known as an Ax-Lindemann-type result; such results are the topic of this section.
\subsection{Ax-Lindemann for \texorpdfstring{$j$}{j}}
In the classical setting, Pila proved the upcoming result in his paper towards Andr\'e-Oort, \cite{Pila2011}.  It is called the Ax-Lindemann theorem for $j$.  Of great interest in its own right, it is also vital to the proof of Andr\'e-Oort-type results, via Pila-Wilkie.  Before we can state it, we will need the following definition:
\begin{definition}
 Consider some subset $Z\suq\uh^n$.  A \emph{complex algebraic component} $A$ of $Z$ is a connected component of a set of the form
 \[W\cap\uh^n,\]
 for $W$ an irreducible subvariety of $\mathbb{C}^n$, with the property that $A\suq Z$.
\end{definition}
\begin{theorem}[Pila, ``Ax-Lindemann for $j$'']\label{thrm:PilaAL}
 Let $V\suq\mathbb{C}^n$ be a variety.  Define a map $\pi:\uh^n\to\mathbb{C}^n$ by
 \[\pi(\tau_1,\dots,\tau_n)=(j(\tau_1),\dots,j(\tau_n)),\]
 and let $\mathcal{Z}=\pi^{-1}(V)$.  
 
 A maximal complex algebraic component of $\mathcal{Z}$ is weakly $\uh$-special. 
\end{theorem}
As we noted in the introduction, this is loosely saying: ``the only complex algebraic relations between coordinates in $\uh$ that induce algebraic relations between their $j$-images in $\mathbb{C}^n$ are the geodesic relations.''

For our purposes, the Ax-Lindemann theorem for $j$ also tells us the following.
\begin{corollary}\label{cor:realAL}
 Let $V$, $\pi$ and $\mathcal{Z}$ be as in \ref{thrm:PilaAL}.
 Then $\mathcal{Z}^\text{alg}$ is simply the union of all positive-dimensional weakly $\uh$-special subvarieties of $\mathcal{Z}$.
\end{corollary}
To go from \ref{thrm:PilaAL} to \ref{cor:realAL} one just uses the holomorphicity of $j$.  A real semialgebraic arc in $\mathcal{Z}$ is contained in a complex algebraic component of $\mathcal{Z}$ by analytic continuation.
\subsection{Quasimodular Ax-Lindemann}
For the QM function $\chi$, a good portion of the work on Ax-Lindemann results is already done for us.  The upcoming result is due to Pila, in \cite{Pila2013}.  To state it, we will need a definition.

\begin{definition}
 Let $\tau_1,\dots,\tau_n$ be elements of some algebraic function field $\mathbb{C}(W)$.  Then $\tau_1,\dots,\tau_n$ are called \emph{geodesically dependent} if either:
 \begin{itemize}
  \item For some $g\in\gl$ and some $i,j$, we have $\tau_i=g\tau_j$ whenever $\tau_i,\tau_j$ take values in $\uh$, \emph{or}
  \item At least one of the $\tau_i$ is constant.
 \end{itemize}
 Otherwise, the $\tau_i$ are called \emph{geodesically independent}.
\end{definition}

\begin{theorem}[Pila, Ax-Lindemann with Derivatives]\label{QmodAxL}

 Suppose that $\mathbb{C}(W)$ is an algebraic function field and that
 \[\tau_1,\dots,\tau_n\in\mathbb{C}(W)\]
 take values in $\uh$ at some $P\in W$, and are geodesically independent.  Then the 3n functions
 \[j(\tau_1),\dots,j(\tau_n),\qquad j'(\tau_1),\dots,j'(\tau_n),\qquad j''(\tau_1),\dots,j''(\tau_n)\]
 (considered as functions on $W$ locally near $P$) are algebraically independent over $\mathbb{C}(W)$.  
\end{theorem}
For our purposes, we need a slightly stronger formulation of this result.

\begin{theorem}[Ax-Lindemann with Derivatives, Stronger Form]\label{thrm:AxLwDerivs}
 Let $F$ be an irreducible polynomial in $3n+1$ variables over $\mathbb{C}$.  Let $A\suq\uh^n$ be a complex algebraic component and let $G$ be the smallest weakly $\uh$-special variety containing $A$.  Suppose that $G$ is a GUT variety and that 
 \[F(\tau_1,j(\tau_1),j'(\tau_1),j''(\tau_1),\dots,j(\tau_n),j'(\tau_n),j''(\tau_n))=0\]
 for all $(\tau_1,\dots,\tau_n)\in A$.  Then in fact this holds for all $(\tau_1,\dots,\tau_n)\in G$.
\end{theorem}
\begin{proof}
  We will work by induction on $n$.  The case $n=1$ is immediate.
  
  By definition, the algebraic component $A$ is a connected component of some variety $W\suq\mathbb{C}^n$.  Treating $\tau_1,\dots,\tau_n$ as the coordinate functions on $W$, the hypotheses of the theorem imply that 
  \[j(\tau_1),j'(\tau_1),j''(\tau_1),\dots,j(\tau_n),j'(\tau_n),j''(\tau_n),\] treated as functions locally near some $P\in A$, are algebraically dependent over $\mathbb{C}(W)$, whence Theorem \ref{QmodAxL} tells us that the $\tau_i$ are geodesically dependent.
  
  By induction, we may assume that no $\tau_i$ is constant on $A$.  Hence there are $1\leq i,j\leq n$ and $g\in\gl$ such that $\tau_i=g\tau_j$ on $A$.  Since this is a symmetric condition, we may assume that $i\ne 1$.  Then without loss of generality, $i=n$.
  
  Since $G$ is a GUT variety, $g$ is upper triangular.  Hence there are algebraic functions $\phi_1,\phi_2,\phi_3$ (induced by the modular polynomials and their derivatives) such that: 
  \begin{equation}\label{eqn:Phi1}j(\tau_i)=\phi_1(j(\tau_j)),\end{equation}
  \begin{equation}\label{eqn:Phi2}j'(\tau_i)=\phi_2(j(\tau_j),j'(\tau_j)),\end{equation}
  \[\text{and}\]
  \begin{equation}\label{eqn:Phi3}j''(\tau_i)=\phi_3(j(\tau_j),j'(\tau_j),j''(\tau_j)).\end{equation}
  Substituting this into $F$ yields
  \begin{multline*}F\bigl[\tau_1,j(\tau_1),j'(\tau_1),j''(\tau_1),\dots,j(\tau_{n-1}),j'(\tau_{n-1}),j''(\tau_{n-1}),\\\phi_1(j(\tau_j)),\phi_2(j(\tau_j),j'(\tau_j)),\phi_3(j(\tau_j),j'(\tau_j),j''(\tau_j))\bigr]=0\end{multline*}
  whenever $(\tau_1,\dots,\tau_{n-1},g\tau_j)\in A$.  We can then rewrite this as
  \[\sigma(\tau_1,j(\tau_1),j'(\tau_1),j''(\tau_1),\dots,j(\tau_{n-1}),j'(\tau_{n-1}),j''(\tau_{n-1}))=0,\]
  for some algebraic function $\sigma$.  This will hold for all $(\tau_1,\dots,\tau_{n-1})\in A'$, where $A'$ is the projection of $A$ onto the first $n-1$ coordinates.
  
  It is possible that $\sigma$ is the zero function.  If so, then working backwards we see that $F$ vanishes whenever (\ref{eqn:Phi1}), (\ref{eqn:Phi2}) and (\ref{eqn:Phi3}) hold.  In particular, $F$ vanishes whenever $\tau_i=g\tau_j$.  Hence it must vanish on $G$, as required.
  
  If $\sigma\ne 0$, we have more work to do.  There is an irreducible polynomial $p_\sigma$ such that
  \[p_\sigma(\sigma(\mathbf{X}),\mathbf{X})=0\]
  for all $\mathbf{X}$.  In particular,
  \begin{equation}\label{eqn:OneFewer}p_\sigma(0,\tau_1,j(\tau_1),j'(\tau_1),j''(\tau_1),\dots,j(\tau_{n-1}),j'(\tau_{n-1}),j''(\tau_{n-1}))=0\end{equation}
  for all $(\tau_1,\dots,\tau_{n-1})\in A'$.  Note that $p_\sigma(0,\mathbf{X})$ is not the zero polynomial.  
  
  We can now appeal to induction to see that (\ref{eqn:OneFewer}) holds for all \[(\tau_1,\dots,\tau_{n-1})\in G',\] where $G'$ is the projection of $G$ onto its first $n-1$ coordinates.  Putting it in different terms: 0 is a root of
  \begin{equation}\label{eqn:RootPoly}p_\sigma(X,\tau_1,j(\tau_1),\dots,j''(\tau_{n-1}))\end{equation}
  whenever $(\tau_1,\dots,\tau_{n-1})\in G'$.  We can choose a point $\mathbf{p}\in A'$, a $G'$-open neighbourhood $V$ of $\mathbf{p}$ and a complex-open neighbourhood $W$ of 0 such that:
  for all $\mathbf{q}\in V$, the only root of (\ref{eqn:RootPoly}) within $W$ is the root 0.  However, $\sigma(\tau_1,j(\tau_1),\dots,j''(\tau_{n-1}))$ is always a root of (\ref{eqn:RootPoly}).  So for all $(\tau_1,\dots,\tau_{n-1})\in V$, we must have
  \[\sigma(\tau_1,j(\tau_1),\dots,j''(\tau_{n-1}))=0.\]
  By analytic continuation, this holds for all $(\tau_1,\dots,\tau_{n-1})\in G'$.  Recalling the definition of $\sigma$, we get that
  \begin{multline*}F\bigl[\tau_1,j(\tau_1),j'(\tau_1),j''(\tau_1),\dots,j(\tau_{n-1}),j'(\tau_{n-1}),j''(\tau_{n-1}),\\\phi_1(j(\tau_j)),\phi_2(j(\tau_j),j'(\tau_j)),\phi_3(j(\tau_j),j'(\tau_j),j''(\tau_j))\bigr]=0\end{multline*}
  whenever $(\tau_1,\dots,\tau_{n-1})\in G'$.  Hence
  \[F(\tau_1,j(\tau_1),\dots,j''(\tau_{n-1}),j(g\tau_j),j'(g\tau_j),j''(g\tau_j))=0\]
  for all $(\tau_1,\dots,\tau_{n-1})\in G'$.  In other words
  \[F(\tau_1,j(\tau_1),\dots,j''(\tau_n))=0\]
  for all $(\tau_1,\dots,\tau_n)\in G$, as required.
 \end{proof}
For our purposes, we need a version of this result that discusses $j$ and $\chi$, rather than the derivatives of $j$, hence the following corollary.  
\begin{corollary}\label{thrm:AxLChiJAndF}
 Let $F$ be an irreducible polynomial in $3n+1$ variables over $\mathbb{C}$.  Let $A\suq\uh^n$ be a complex algebraic component and let $G$ be the smallest weakly $\uh$-special variety containing $A$.  Suppose that $G$ is a GUT variety and that 
 \[F(\tau_1,j(\tau_1),\chi(\tau_1),f(\tau_1),\dots,j(\tau_n),\chi(\tau_n),f(\tau_n))=0\]
 for all $(\tau_1,\dots,\tau_n)\in A$.  Then in fact this holds for all $(\tau_1,\dots,\tau_n)\in G$.
 (Recall that $f$ is the function $E_4E_6/\Delta$, which arises in the transformation law for $\chi$ and as the coefficient of $1/\Imm\tau$ in $\chi^*$.)
 \end{corollary}
 \begin{proof}
  Follows easily from \ref{thrm:AxLwDerivs}, using the fact that $j,\chi,f\in\mathbb{C}(j,j',j'')$ and that $j(\tau),\chi(\tau),f(\tau)$ are algebraically independent functions over $\mathbb{C}(\tau)$.
 \end{proof}

\subsection{Almost Holomorphic Ax-Lindemann}
In the classical situation, as we see above, the holomorphicity of the functions involved allows us to `complexify the parameter' to produce a complex algebraic set from a real algebraic one.  Since $\chi^*$ is not holomorphic, there is substantial difficulty in attempting to complexify the parameter in the same way.  While a real algebraic arc in $\uh^n$ is certainly contained in a complex algebraic component of $\uh^n$, there is no guarantee that this algebraic component remains within the preimage of the given variety $V$.  Fortunately, the simple shape of $\chi^*$ allows us to use some tricks to get around this problem.  This subsection is dedicated to proving the desired Ax-Lindemann results for $\chi^*$.  This is a crucial step towards our central Andr\'e-Oort result for $\chi^*$; most of the novelty in our proof of  \ref{thrm:AOforChiStar2} lies in this nonholomorphic Ax-Lindemann result.

As we mentioned in section \ref{sect:intro}, we will be discussing a map $\pi:\uh^n\to\mathbb{C}^{2n}$, defined by
\[\pi(\tau_1,\dots,\tau_n)=(j(\tau_1),\chi^*(\tau_1),\dots,j(\tau_n),\chi^*(\tau_n)).\]

\begin{theorem}[AHM Ax-Lindemann]\label{thrm:AxLChiStarAndJ}
 Let $S$ be an arc of a real algebraic curve in $\uh^n$ and suppose that $S\suq \pi^{-1}(V)$, where $V$ is some irreducible variety in $\mathbb{C}^{2n}$.  Then $S$ is contained in a weakly $\uh$-special variety $G$ with $G\suq \pi^{-1}(V)$.
\end{theorem}
The proof of this is necessarily rather technical, so for ease of reading we have broken it into various smaller chunks.  The plan is as follows.  Firstly, we deal with the case in which the imaginary part of every complex coordinate is constant on the arc $S$.  This is the content of Lemma \ref{lma:RemovingFs} and Corollary \ref{lma:ConstantIms}.  With this done, we can assume that the imaginary part of at least one coordinate (say $\tau_1$) is nonconstant on $S$.  Hence we can parametrise $S$ in terms of the imaginary part of $\tau_1$.  

Using this parametrisation, we will show that a particular algebraic function $\phi$ in the variables $\Imm\tau_1$, $j(\tau_i)$, $\chi(\tau_i)$, $f(\tau_i)$ vanishes on $S$.  If $\phi$ takes a very specific shape, we can conclude via \ref{thrm:AxLChiJAndF}.  Otherwise, we will see that
\[\Imm\tau_1=\psi(j(\tau_1),\dots,j(\tau_n),\chi(\tau_1),\dots,\chi(\tau_n),f(\tau_1),\dots,f(\tau_n))\]
on $S$, for some algebraic function $\psi$.  In this situation, Lemma \ref{lma:yEqualsAlgebraic} shows that $\Imm\tau_1$ must be constant on $S$ after all, which is a contradiction.

\begin{lemma}\label{lma:RemovingFs}
 Let $G\suq\uh^n$ be a GUT variety, let $F$ be a polynomial in $2n$ variables, and let $c_1,\dots,c_n$ be real constants.  Suppose that
 \[F\left(j(\tau_1),\chi(\tau_1)-\dfrac{3}{\pi c_1}f(\tau_1),\dots,j(\tau_n),\chi(\tau_n)-\dfrac{3}{\pi c_n}f(\tau_n)\right)=0\]
 for all $(\tau_1,\dots,\tau_n)\in G$.  Then 
 \[F\left(j(\tau_1),\chi^*(\tau_1),\dots,j(\tau_n),\chi^*(\tau_n)\right)=0\]
 for all $(\tau_1,\dots,\tau_n)\in G$.
 \end{lemma}
\begin{proof}
  By induction on $n$ we may assume that no coordinate is constant on $G$.  So up to permutation of coordinates, we have
  \begin{equation*}G=\bigl\{(\tau_1,g_{1,1}\tau_1,\dots,g_{1,k_1}\tau_1,\tau_2,g_{2,1}\tau_2,\dots,g_{2,k_2}\tau_2,\dots,\tau_r,g_{r,1}\tau_r,\dots,g_{r,k_r}\tau_r):\tau_1,\dots,\tau_r\in\uh\bigr\},\end{equation*}
  for some upper triangular matrices $g_{i,j}$.
  Hence
  \begin{multline}\label{eqn:FZeroOnG}
   F\bigg[j(\tau_1),\chi(\tau_1)-\dfrac{3}{\pi d_1}f(\tau_1),\dots,j(g_{1,k_1}\tau_1),\chi(g_{1,k_1}\tau_1)-\dfrac{3}{\pi d_{1,k_1}}f(g_{1,k_1}\tau_1),\qquad\\
   \dots,\\
   j(\tau_r),\chi(\tau_r)-\dfrac{3}{\pi d_r}f(\tau_r),\dots,j(g_{r,k_r}\tau_r),\chi(g_{r,k_r}\tau_r)-\dfrac{3}{\pi d_{r,k_r}}f(g_{r,k_r}\tau_r)\bigg]=0,  
  \end{multline} 
  for some suitable relabelling $d_i$, $d_{i,j}$ of the constants $c_i$.
  
  All of the $g_{i,j}$ are upper triangular matrices in $\gl$, so let us consider a general upper triangular matrix $g=\begin{pmatrix}a&b\\0&d\end{pmatrix}$.  Let $A=\gcd(b,d)$ and $D=ad/A$.  Let $k,m$ be integers such that $mb+kd=A$.  For all integers $t$, we have
  \[\begin{pmatrix}b/A&-k+tb\\d/A&m+td\end{pmatrix}\cdot\begin{pmatrix}A&-ma\\0&D\end{pmatrix}=\begin{pmatrix}a&b\\0&d\end{pmatrix}\cdot\begin{pmatrix}0&-1\\1&tD\end{pmatrix}.\]
  The leftmost matrix is an element of $\slz$.  The matrix $\begin{pmatrix}A&-ma\\0&D\end{pmatrix}$ has the same determinant as $g$; we shall call this matrix $h$.  Note (taking $t=0$ above) that
  \[\begin{pmatrix}b/A&-k\\d/A&m\end{pmatrix}h\begin{pmatrix}0&1\\-1&0\end{pmatrix}=g.\]
  From these matrix equations and the transformation properties of $j$, $\chi$ and $f$, we can easily see that
  \begin{align*}j\left(g\begin{pmatrix}0&-1\\1&tD\end{pmatrix}\tau\right)&=j(h\tau), &\text{ for all }t.\\
  \chi\left(g\begin{pmatrix}0&-1\\1&tD\end{pmatrix}\tau\right)&\to\chi(h\tau)&\text{ as }t\to\infty.\\
  f\left(g\begin{pmatrix}0&-1\\1&tD\end{pmatrix}\tau\right)&\to 0&\text{ as }t \to \infty.
  \end{align*}  
  Also,
  \begin{align*}j\left(\begin{pmatrix}0&-1\\1&tD\end{pmatrix}\tau\right)&=j(\tau).&\\
  \chi\left(\begin{pmatrix}0&-1\\1&tD\end{pmatrix}\tau\right)&\to \chi(\tau)&\text{ as }t\to\infty.\\f\left(\begin{pmatrix}0&-1\\1&tD\end{pmatrix}\tau\right)&\to 0&\text{ as }t\to \infty.\end{align*}
  Now, equation (\ref{eqn:FZeroOnG}) holds for all $\tau_1,\dots,\tau_r\in\uh$.  Hence we can replace each $\tau_i$ in (\ref{eqn:FZeroOnG}) by $\begin{pmatrix}0&-1\\1&tD_i\end{pmatrix}\tau_i$, for suitable fixed $D_i$ and arbitrary $t$.  Letting $t$ tend to infinity we see by continuity of $F$ that
  \begin{equation*}
   F\bigg[j(\tau_1),\chi(\tau_1),\dots,j(h_{1,k_1}\tau_1),\chi(h_{1,k_1}\tau_1),
   \dots,
   j(\tau_r),\chi(\tau_r),\dots,j(h_{r,k_r}\tau_r),\chi(h_{r,k_r}\tau_r)\bigg]=0 
  \end{equation*}
  for all $\tau_i\in\uh$ and certain upper triangular matrices $h_{j,k}$.  By \ref{lma:isomorphism} (the isomorphism between upper triangular extensions of the fields of QM/AHM functions), we therefore have
  \begin{equation}\label{eqn:FWithChiStar}
   F\bigg[j(\tau_1),\chi^*(\tau_1),\dots,j(h_{1,k_1}\tau_1),\chi^*(h_{1,k_1}\tau_1),
   \dots,
   j(\tau_r),\chi^*(\tau_r),\dots,j(h_{r,k_r}\tau_r),\chi^*(h_{r,k_r}\tau_r)\bigg]=0
  \end{equation}  
  The matrices $h_{j,k}$ each have the same relation to $g_{j,k}$ as $h$ does to $g$ in the calculation above.  In particular, there is $\gamma_{j,k}\in\slz$ such that
  \[\gamma_{j,k}h_{j,k}\begin{pmatrix}0&1\\-1&0\end{pmatrix}=g_{j,k}.\]
  So we can replace each $\tau_i$ in (\ref{eqn:FWithChiStar}) by $\begin{pmatrix}0&1\\-1&0\end{pmatrix}\tau_i$ and use the modularity of $j$ and $\chi^*$ to see that
  \begin{equation*}
   F\bigg[j(\tau_1),\chi^*(\tau_1),\dots,j(g_{1,k_1}\tau_1),\chi^*(g_{1,k_1}\tau_1),
   \dots,
   j(\tau_r),\chi^*(\tau_r),\dots,j(g_{r,k_r}\tau_r),\chi^*(g_{r,k_r}\tau_r)\bigg]=0.
  \end{equation*} 
 This says precisely that \[F(j(\tau_1),\chi^*(\tau_1),\dots,j(\tau_n),\chi^*(\tau_n))=0\] for all $(\tau_1,\dots,\tau_n)\in G$.
 \end{proof}
 
\begin{corollary}\label{lma:ConstantIms}
 Let $S$ and $V$ be as in \ref{thrm:AxLChiStarAndJ}.  Suppose that the imaginary part of every complex coordinate is constant on $S$.  Then $S$ is contained in a weakly $\uh$-special variety $G$ with $G\suq \pi^{-1}(V)$.
\end{corollary}
\begin{proof}  
  By induction on $n$, we may assume that no complex coordinate is constant on $S$.  So consider the smallest weakly $\uh$-special variety containing $S$, which we will call $G$.  Since no complex coordinate is constant on $S$, the same is true of $G$.  We want to show $G\suq\pi^{-1}(V)$.
  
  Consider some coordinate $\tau_i$ on $S\suq G$.  It takes the form $\tau_i=x_i+ic_i$.  Suppose that on $G$, there is some $\tau_j$, $j\ne i$ which is related to $\tau_i$ by some matrix $g$ which fails to be upper triangular.  Then on $S$, we have $\tau_j=g(x_i+ic_i)$.  Since $\tau_i$ is nonconstant on $S$, $x_i$ must vary, which then forces $\Imm\tau_j$ to vary since $g$ is not upper triangular.  This is a contradiction.  So $G$ is a GUT variety.  
  
  Now pick any of the irreducible polynomials $F$ which define $V$.  We have
  \begin{equation}\label{eqn:FVanishesWithConsts}F\left(j(\tau_1),\chi(\tau_1)-\dfrac{3}{\pi c_1}f(\tau_1),\dots,j(\tau_n),\chi(\tau_n)-\dfrac{3}{\pi c_n}f(\tau_n)\right)=0\end{equation}
  for all $(\tau_1,\dots,\tau_n)\in S$ and for real constants $c_i=\Imm\tau_i$.  
  
  Let us parametrise $S$ in terms of some real parameter $t$, as the image of a map $t\mapsto(\tau_1(t),\dots,\tau_n(t))$ around $t=0$.  Without loss of generality, suppose that $\tau_1$ is nonconstant, so that all of the other functions $\tau_i$ are algebraic over $\tau_1$.  The functions $\tau_i$ may then be extended to complex $t$ in some complex neighbourhood of $0$.  The image of this complex neighbourhood under the map then necessarily lives in some irreducible complex algebraic curve $C$.  Since (\ref{eqn:FVanishesWithConsts}) holds on $S\suq C$ and all of the functions arising in (\ref{eqn:FVanishesWithConsts}) are complex analytic, it follows that (\ref{eqn:FVanishesWithConsts}) holds on the whole of $C$.  (This method of complexifying the parameter will arise several times; compare with, for instance, \cite[Lemma 2.1]{Pila2009}.)
  
  So we get that (\ref{eqn:FVanishesWithConsts}) holds on some complex algebraic component $A$ containing $S$.  Define $G$ to be the smallest weakly $\uh$-special variety containing $A$.  As previously, we may assume that $G$ is a GUT variety.  Hence we can apply \ref{thrm:AxLChiJAndF} to see that 
  \[F\left(j(\tau_1),\chi(\tau_1)-\dfrac{3}{\pi c_1}f(\tau_1),\dots,j(\tau_n),\chi(\tau_n)-\dfrac{3}{\pi c_n}f(\tau_n)\right)=0\]
  for all $(\tau_1,\dots,\tau_n)\in G$.  By Lemma \ref{lma:RemovingFs}, we then have
  \[F\left(j(\sigma_1),\chi^*(\sigma_1),\dots,j(\sigma_n),\chi^*(\sigma_n)\right)=0\]
  for all $(\sigma_1,\dots,\sigma_n)\in G$.  This holds for all of the defining polynomials of $V$, hence $S\suq G\suq\pi^{-1}(V)$ as required.
 \end{proof}

\begin{lemma}\label{lma:yEqualsAlgebraic}
 Let $S$ be an arc of a real algebraic curve in $\uh^n$ and let $\psi$ be an algebraic function in $3n$ variables.  Suppose that
 \[\Imm\tau_1=\psi(j(\tau_1),\dots,j(\tau_n),\chi(\tau_1),\dots,\chi(\tau_n),f(\tau_1),\dots f(\tau_n))\]
 for all $(\tau_1,\dots,\tau_n)\in S$.  
 
 Let $G$ be the smallest weakly $\uh$-special variety containing $S$, and suppose that $G$ is a GUT variety.  Then $\Imm\tau_1$ is constant on $S$.
\end{lemma}
 \noindent\textbf{Notation:} The tuple \[(j(\tau_1),\dots,j(\tau_n),\chi(\tau_1),\dots,\chi(\tau_n),f(\tau_1),\dots f(\tau_n))\] will arise often in what follows, so we abbreviate it as $\tilde{\pi}(\tau_1,\dots,\tau_n)$.  We will also abbreviate $y=\Imm\tau_1$ throughout.  So the first hypothesis of the Lemma may be written as
  \[y=\psi(\tilde\pi(\tau_1,\dots,\tau_n)).\] 
 \begin{proof}[Proof of \ref{lma:yEqualsAlgebraic}]
  Suppose for a contradiction that $y=\Imm\tau_1$ is nonconstant on $S$.  Then we can parametrise $S$ in terms of $y$, yielding
  \[S=\{(x(y)+iy,u_2(y)+iv_2(y),\dots,u_n(y)+iv_n(y)):y\in U\}\]
  for some open set $U\suq\mathbb{R}$ and algebraic functions $x,u_i,v_i$, real-valued on $U$.
  
  Since $S$ is an algebraic arc, we also have some polynomials $a_i$ such that
 \[a_i(x(y),y,u_2(y),v_2(y),\dots,u_n(y),v_n(y))=0\]
 for all $(\tau_1,\dots,\tau_n)\in S$.  Noting that $\tau_1=x(y)+iy$, and replacing instances of $y$ with $\psi$, we get
 \begin{multline}\label{eqn:aIEqualsZero}
  a_i\bigl[\tau_1-i\psi(\tilde\pi(\tau_1,\dots,\tau_n)),
      \psi(\tilde\pi(\tau_1,\dots,\tau_n)),\\
      u_2(\psi(\tilde\pi(\tau_1,\dots,\tau_n))),
      v_2(\psi(\tilde\pi(\tau_1,\dots,\tau_n))),\\
      \dots,\\
      u_n(\psi(\tilde\pi(\tau_1,\dots,\tau_n))),v_n(\psi(\tilde\pi(\tau_1,\dots,\tau_n)))\bigr]=0,
 \end{multline}
 for all $(\tau_1,\dots,\tau_n)\in S$.
 We rewrite the left hand side of this equation as an algebraic function 
 \[\sigma(\tau_1,\tilde\pi(\tau_1,\dots,\tau_n)).\]
 Then there is an irreducible polynomial $p_\sigma$ such that
 \[p_\sigma(\sigma(T,\mathbf{J},\mathbf{X},\mathbf{F}),T,\mathbf{J},\mathbf{X},\mathbf{F})=0\]
 identically.  In particular, since $\sigma$ vanishes on $S$, we have
 \[p_\sigma(0,\tau_1,\tilde\pi(\tau_1,\dots,\tau_n))=0\]
 for all $(\tau_1,\dots,\tau_n)\in S$.  By complexifying the parameter, as in Corollary \ref{lma:ConstantIms} and \cite[Lemma 2.1]{Pila2009}, this holds on a complex algebraic component $A$ containing $S$.  Now, the weakly special closure of $A$ is the same as the weakly special closure of $S$, namely $G$.  Since $G$ is a GUT variety, we may therefore apply \ref{thrm:AxLChiJAndF} to see that 
 \[p_\sigma(0,\tau_1,\tilde\pi(\tau_1,\dots,\tau_n))=0\]
 for all $(\tau_1,\dots,\tau_n)\in G$.
 
 In other words, 0 is a root of
 \begin{equation}\label{eqn:pSigma}p_\sigma(X,\tau_1,\tilde\pi(\tau_1,\dots,\tau_n))\end{equation}
 for all $(\tau_1,\dots,\tau_n)\in G$.  Since
 \[p_\sigma(\sigma(T,\mathbf{J},\mathbf{X},\mathbf{F}),T,\mathbf{J},\mathbf{X},\mathbf{F})=0\]
 identically, we know that 
 \[\sigma(\tau_1,\tilde\pi(\tau_1,\dots,\tau_n))\]
 is also root of (\ref{eqn:pSigma}) for all $(\tau_1,\dots,\tau_n)\in G$.  
 
 We can pick a point $a\in S$, a $G$-open neighbourhood $W$ of $a$, and a complex neighbourhood $U$ of 0, such that: as $(\tau_1,\dots,\tau_n)$ varies within $W$, the only root of (\ref{eqn:pSigma}) within $U$ is 0 itself.  However, as $(\tau_1,\dots,\tau_n)$ varies in $W$, the function $\sigma(\tau_1,\tilde{\pi}(\tau_1,\dots,\tau_n))$ remains a root of (\ref{eqn:pSigma}).  Since $\sigma$ vanishes on $S$, we can get it arbitrarily close to 0 within $W$.  In particular, we can get $\sigma(\tau_1,\tilde{\pi}(\tau_1,\dots,\tau_n))$ to lie within $U$.  Since it is a root of (\ref{eqn:pSigma}), we must have 
 \[\sigma(\tau_1,\tilde{\pi}(\tau_1,\dots,\tau_n))=0\]
 for all $(\tau_1,\dots,\tau_n)\in W$.  By analytic continuation, this holds everywhere on $G$, which says that (\ref{eqn:aIEqualsZero}) holds on $G$.
 
 For notational simplicity, let us suppose that the coordinates which are related to $\tau_1$ in $G$ are the first $k$ coordinates, that is:
 \[G=\{(\tau_1,g_2\tau_1,\dots,g_k\tau_1):\tau_1\in\mathcal{H}\}\times G',\]
 for some other GUT variety $G'$.
 So, whenever \[(\tau_1,g_2\tau_1,\dots,g_k\tau_1,\tau_{k+1},\dots,\tau_n)\in G,\] we also have
 \[(\tau_1+t,g_2(\tau_1+t),\dots,g_k(\tau_1+t),\tau_{k+1},\dots\tau_n)\in G\]
 for every $t\in \mathbb{Z}$.  Since $G$ is a GUT set, the $g_i$ are upper triangular, so the numbers $g_i(\tau_1+t)$, up to translation by an integer, take only finitely many values as $t$ varies.  In particular, since $j$, $\chi$ and $f$ are periodic, each of the functions
 \[j(g_i(\tau_1+t)),\qquad\chi(g_i(\tau_1+t)),\qquad f(g_i(\tau_1+t))\]
 takes only finitely many values as $t$ varies.  Hence
 \begin{align*}
  \psi\big(\tilde\pi(\tau_1+t,g_2(\tau_1+t),\dots,g_k(\tau_1+t),\tau_{k+1},\dots,\tau_n)\big)
 \end{align*}
 takes only finitely many values as $t\in\mathbb{Z}$ varies.  If we plug this into (\ref{eqn:aIEqualsZero}), we see that 
 \[a_i(\tau_1+t-ic,c,u_2(c),v_2(c),\dots,v_n(c))=0\]
 for some constant $c$ and infinitely many distinct $t$.  Thus $a_i$ is independent of its first coordinate.  Since this is true of all the $a_i$ defining $S$, the only possibility for $S$ is that it is the product of a horizontal line in the $\tau_1$ plane and points in the other coordinates.  So $y$ is constant on $S$, which is a contradiction.
 \end{proof}
With all the above lemmas done, we may finally proceed to the body of the proof of \ref{thrm:AxLChiStarAndJ}.

\begin{proof}[\textbf{Proof of \ref{thrm:AxLChiStarAndJ}}]
 By induction on $n$, we may assume that no complex coordinate is constant on $S$.  It might be, however, that the imaginary part of one or more coordinates is constant on $S$.  If $\Imm\tau_i$ is constant on $S$ for every $i$, then we are in the situation of Lemma \ref{lma:ConstantIms} so we conclude immediately.  Hence we may assume without loss of generality that $\Imm\tau_1$ is nonconstant on $S$.  
 
 Next, let $G$ be the unique smallest weakly $\uh$-special subvariety of $\mathcal{H}^n$ containing $S$.  It is a standard fact (which we have used once already; see for instance Lang \cite{Lang1987} or \cite[Exercise 1.2.11]{Diamond2005}) that any $g\in\gl$ takes the form $\gamma\cdot h$ for some upper triangular $h\in\gl$ and some $\gamma\in\slz$.  Therefore there is some $\gamma\in\slz^n$ such that $\gamma G$ is a GUT variety.  The subset $\gamma S\suq \gamma G$ is still a real semialgebraic arc.  By the modularity of $j$ and $\chi^*$, $\gamma G\suq\pi^{-1}(V)$ if and only if $G\suq\pi^{-1}(V)$.  So by working with $\gamma S$ we may assume without loss of generality that $G$ is a GUT variety.
 
 We will write $y=\Imm\tau_1$ throughout, and retain the abbreviation
 \[\tilde\pi(\tau_1,\dots,\tau_n)=(j(\tau_1),\dots,j(\tau_n),\chi(\tau_1),\dots,\chi(\tau_n),f(\tau_1),\dots,f(\tau_n)).\] 
 Since $y$ is nonconstant, we can parametrise $S$ as
 \[S=\{(x(y)+iy,u_2(y)+iv_2(y),\dots,u_n(y)+iv_n(y)):y\in U\},\]
 for some open $U\suq \mathbb{R}$ and algebraic functions $x$, $u_i$, $v_i$, real-valued on $U$.
 
 Consider one of the polynomials $F$ which defines $V$.  We have
 \[F\left(j(\tau_1),\chi(\tau_1)-\dfrac{3}{\pi y}f(\tau_1),\dots,j(\tau_n),\chi(\tau_n)-\dfrac{3}{\pi v_n(y)}f(\tau_n)\right)=0\]
 for all $(\tau_1,\dots,\tau_n)\in S$.  We can rewrite the left hand side of this equation as an algebraic function
 \[\phi(y,\tilde\pi(\tau_1,\dots,\tau_n)).\] 
 Since $\phi$ is an algebraic function, there is an irreducible polynomial $p_\phi$ with the property that
 \[p_\phi(\phi(T,\mathbf{J},\mathbf{X},\mathbf{F}),T,\mathbf{J},\mathbf{X},\mathbf{F})=0\]
 for all $T$, $\mathbf{J}=(J_1,\dots,J_n)$, $\mathbf{X}=(X_1,\dots,X_n)$ and $\mathbf{F}=(F_1,\dots,F_n)$.  In particular, we have that
 \[p_\phi(0,y,\tilde{\pi}(\tau_1,\dots,\tau_n))=0\]
 for all $(\tau_1,\dots,\tau_n)\in S$.  So let us define
 \[P(T,\mathbf{J},\mathbf{X},\mathbf{F})=p_\phi(0,T,\mathbf{J},\mathbf{X},\mathbf{F}).\]
 Note that $P$ is not the zero polynomial, since $p_\phi$ is irreducible.
 
 We are going to modify $P$ as follows.  Consider each coefficient of $T^k$ in $P$ separately.  These are polynomials 
 \[C_k(\mathbf{J},\mathbf{X},\mathbf{F}).\]
 For each $k$, if 
 \[C_k(\tilde{\pi}(\tau_1,\dots,\tau_n))=0\]
 for $(\tau_1,\dots,\tau_n)\in S$, then remove this coefficient of $T^k$ from the polynomial $P$.  Having done this for each coefficient, we have a modified polynomial which we call $\tilde P$.  Note that we still have
 \[\tilde P(y,\tilde{\pi}(\tau_1,\dots,\tau_n))=0\]
 for $(\tau_1,\dots,\tau_n)\in S$.
 
 It is possible that $\tilde P$ is the zero polynomial.  This happens if and only if every coefficient $C_k$ has the property that
 \begin{equation}\label{eqn:CoeffVanish}C_k(\tilde{\pi}(\tau_1,\dots,\tau_n))=0\end{equation}
 for $(\tau_1,\dots,\tau_n)\in S$.  By complexifying the parameter, as in Corollary \ref{lma:ConstantIms} and \cite[Lemma 2.1]{Pila2009}, the equation (\ref{eqn:CoeffVanish}) holds for $(\tau_1,\dots,\tau_n)\in A$, where $A$ is the smallest complex algebraic component containing $S$.  Now, the weakly $\uh$-special closure of $A$ is the same as the weakly special closure of $S$, which is the GUT variety $G$.  Hence we can apply \ref{thrm:AxLChiJAndF}, to see that (\ref{eqn:CoeffVanish}) holds for all $(\tau_1,\dots,\tau_n)\in G$.  
 
 Since this holds for all $C_k$, we have
 \begin{align*}
  &p_\phi(0,Y,\tilde{\pi}(\tau_1,\dots,\tau_n))\\
  =\:&P(Y,\tilde{\pi}(\tau_1,\dots,\tau_n))\\=\:&0
 \end{align*}
 for all $(\tau_1,\dots,\tau_n)\in G$ and all choices of $Y$.  In other words, 0 is a root of
 \begin{equation}\label{eqn:pPhiIsZero}p_\phi(X,Y,\tilde{\pi}(\tau_1,\dots,\tau_n))\end{equation}
 for all $(\tau_1,\dots,\tau_n)\in G$ and all $Y$.
 
 Now we proceed exactly as we did in Lemma \ref{lma:yEqualsAlgebraic}.  We can certainly pick a point $a=(a_1,\dots,a_n)\in S$ such that for all $(\tau_1,\dots,\tau_n)$ in some $G$-open neighbourhood $W$ of $a$, the only root of (\ref{eqn:pPhiIsZero}), in some complex neighbourhood $U$ of 0, is 0 itself.  However, we know that 
 \[X=\phi(Y,\tilde{\pi}(\tau_1,\dots,\tau_n))\]
 is a root of (\ref{eqn:pPhiIsZero}) identically.  Fixing $Y=\Imm a_1$, we see that \[\phi(\Imm a_1,\tilde\pi(\tau_1,\dots,\tau_n))\] gets arbitrarily close to 0 within $W$ (it vanishes at $a$).  So as $(\tau_1,\dots,\tau_n)$ varies within $W$, $\phi$ is a root of (\ref{eqn:pPhiIsZero}), and lies inside of $U$.  The only such root is 0, so we must have
 \[\phi(\Imm a_1,\tilde\pi(\tau_1,\dots\tau_n))=0\]
 for all $(\tau_1,\dots,\tau_n)\in W$.  By analytic continuation, this holds for all \[(\tau_1,\dots,\tau_n)\in G.\]  Recalling the definition of $\phi$, we get
 \[F\left(j(\tau_1),\chi(\tau_1)-\dfrac{3}{\pi \Imm a_1}f(\tau_1),\dots,j(\tau_n),\chi(\tau_n)-\dfrac{3}{\pi v_n(\Imm a_1)}f(\tau_n)\right)=0\]
 for all $(\tau_1,\dots,\tau_n)\in G$.  Hence we are in the situation of Lemma \ref{lma:RemovingFs}, so we get
 \[F\left(j(\tau_1),\chi^*(\tau_1),\dots,j(\tau_n),\chi^*(\tau_n)\right)=0,\]
 for all $(\tau_1,\dots,\tau_n)\in G$, as required.
 
 \bigskip\noindent
 We have now dealt with the case where $\tilde P$ is the zero polynomial.  So we suppose that $\tilde P\ne 0$ and look for a contradiction.  Since
 \[\tilde P(y,\tilde{\pi}(\tau_1,\dots,\tau_n))=0\]
 for $(\tau_1,\dots,\tau_n)\in S$, there is an irreducible factor $Q$ of $\tilde P$ with this same property.
 
 Suppose some coefficient of $y^k$ in $Q$ vanishes on $S$.  Then we repeat the entire process, removing redundant coefficients to get a polynomial $\tilde Q$.  Again, an irreducible component of $\tilde{Q}$ must vanish on $S$.  Then we can remove redundant coefficients from this irreducible component, and so on.
 
 We continue repeating this process until it terminates with an irreducible polynomial
 \[R(y,\tilde{\pi}(\tau_1,\dots,\tau_n)),\]
 which vanishes on $S$, with the property that none of the coefficients of $y^k$ in $R$ vanish on $S$.  If $R$ were the zero polynomial, then working backwards we see that $\tilde P$ should have been the zero polynomial, which we have assumed is not the case.  So $R\ne 0$.  In particular, $R$ is nonconstant as a polynomial in $y$.  
  
 Hence, since none of the coefficients of $y^k$ in $R$ vanish on $S$, we can extract an algebraic function $\psi$ such that
 \[y=\psi(\tilde{\pi}(\tau_1,\dots,\tau_n))\]
 for all $(\tau_1,\dots,\tau_n)\in S$.  By our earlier comment, we know that the smallest $\uh$-special variety containing $S$ is $G$, a GUT variety.  So we are in the situation of Lemma \ref{lma:yEqualsAlgebraic}, hence $y$ is constant on $S$, which is a contradiction. 
\end{proof}

We can reformulate \ref{thrm:AxLChiStarAndJ} into the following slightly cleaner statement.
\begin{corollary}\label{cor:ChiStarALReformed}
 Let $V$ be an irreducible subvariety of $\mathbb{C}^{2n}$ and let $\mathcal{Z}=\pi^{-1}(V)$.  Then $\mathcal{Z}^\text{alg}$ is just the union of the weakly $\uh$-special subvarieties of $\mathcal{Z}$.
\end{corollary}
In the next section we use this to prove the central result of the document.
\section{Andr\'e-Oort for \texorpdfstring{$\chi^*$}{Chi-Star}}\label{sect:AO}
Since there are no obvious $\chi$- or $(j,\chi)$-special varieties in $\mathbb{C}^n$, it is not clear what an Andr\'e-Oort statement should look like.  We can, however, formulate meaningful Andr\'e-Oort statements for $\chi^*$.  In this section we state and prove \ref{thrm:AOforChiStar2}, which is the main theorem of the document, an Andr\'e-Oort theorem for $j$ and $\chi^*$.  The map $\pi$ will throughout be defined as before, namely
\[\pi(\tau_1,\dots,\tau_n)=(j(\tau_1),\chi^*(\tau_1),\dots,j(\tau_n),\chi^*(\tau_n)).\]

The proof follows the standard strategy explicated in \cite{Pila2011} very closely, and we will borrow ideas freely from there.  Readers familiar with the strategy will be aware of the piece that is currently missing.  We need some number-theoretic lower bound in order to force $\pi^{-1}(V)$ to contain many points of a given height.  This will force a real algebraic arc to exist in $\pi^{-1}(V)$, so that we can apply the results of the previous section.  For the case of $j$, the lower bound comes from the size of certain Galois orbits, which are known by a result of Siegel to be sufficiently large.  Our approach essentially comes down to that same lower bound of Siegel, but first we have to do some work to ensure that the bound still applies to $\chi^*$-special points.

\begin{proposition}[Masser]
 For a quadratic point $\tau\in\uh$, we have \[\mathbb{Q}(\chi^*(\tau))\suq\mathbb{Q}(j(\tau)).\]
\end{proposition}
\begin{proof}
  Masser proves this in the Appendix of \cite{Masser1975} for a function he calls $\psi$, which is $E_2^*E_4/E_6$.  Since $\chi^*$ lies in $\mathbb{Q}(\psi, j)$, the result follows for $\chi^*$.
\end{proof}

A careful look at Masser's proof of the above yields the following stronger result.

\begin{proposition}\label{propn:GaloisOrbits}
 Let $\tau\in\uh$ be a quadratic point and consider the algebraic numbers $j(\tau)$ and $\chi^*(\tau)$.  Let $\sigma$ be a Galois conjugation acting on $\mathbb{Q}(j(\tau))\supseteq\mathbb{Q}(\chi^*(\tau))$.  Let $\tau'$ be a quadratic point such that $j(\tau')=\sigma(j(\tau))$.  Then $\chi^*(\tau')=\sigma(\chi^*(\tau))$.
\end{proposition}
\begin{proof}
  This comes entirely from close inspection of Masser's work (the appendix in \cite{Masser1975}).  
  Let $d$ be the discriminant of the quadratic number $\tau$, and suppose that $d$ is not equal to $3k^2$ for some odd $k$.  Define some rational functions $\beta_{i,k}^\tau$ such that $\beta_{i,k}^{\tau}(j(\tau))$ are the coefficients of the Taylor expansion of $\Phi_d$ about the point $(j(\tau),j(\tau))$.  This we can certainly do, and we get
  \[\Phi_d(X,Y)=\sum_{(i,k)\ne(0,0)}\beta_{i,k}^{\tau}(j(\tau))(X-j(\tau))^i(Y-j(\tau))^k.\]
  It appears that the rational functions $\beta_{i,k}^{\tau}$ will differ with $\tau$.  However, we will show that, for the $\tau$ and $\tau'$ defined in the hypotheses of the theorem, we do have $\beta_{i,k}^{\tau}=\beta_{i,k}^{\tau'}$.
  
  Since $\Phi_d$ has rational coefficients, any Galois conjugation preserves the left hand side of the above.  So we get
  \begin{align*}
   \Phi_d(X,Y)&=\sum\sigma(\beta_{i,k}^{\tau}(j(\tau)))(X-\sigma(j(\tau)))^i(Y-\sigma(j(\tau)))^k\\
   &=\sum\beta_{i,k}^{\tau}(j(\tau'))(X-j(\tau'))^i(Y-j(\tau'))^k.
  \end{align*}
  We also have 
  \[\Phi_d(X,Y)=\sum_{(i,k)\ne(0,0)}\beta_{i,k}^{\tau'}(j(\tau'))(X-j(\tau'))^i(Y-j(\tau'))^j,\]
  so by uniqueness of Taylor coefficients, the rational functions $\beta_i^{\tau}$ and $\beta_j^{\tau'}$ are equal.
  On pages 118 and 119 of \cite{Masser1975}, $\psi(\tau)$ is expressed as a fixed $\mathbb{Q}$-rational function $p$ in the $\beta_{i,k}^{\tau}(j(\tau))$ and $j(\tau)$.  The equality
  \[\psi(\tau)=p(j(\tau),\beta_{i,k}^{\tau}(j(\tau)))\]
  holds whenever $\tau$ has discriminant $d$ and $\beta_{i,j}^{\tau}$ are the Taylor coefficients of $\Phi_d$ about $(j(\tau),j(\tau))$.  Since $\tau'$ and $\tau$ have the same discriminant (both satisfy $\Phi_d(j(\rho),j(\rho))=0$), this equation holds for both $\tau$ and $\tau'$.  Since $\beta_i^\tau=\beta_i^{\tau'}$ we get
  \begin{align*}\sigma\psi(\tau)&=p(\sigma(j(\tau)),\beta_{i,k}^{\tau}(\sigma(j(\tau))))\\
   &=p(j(\tau'),\beta_{i,k}^{\tau}(j(\tau')))\\
   &=p(j(\tau'),\beta_{i,k}^{\tau'}(j(\tau')))\text{ since }\beta_{i,k}^{\tau}=\beta_{i,k}^{\tau'}\\
   &=\psi(\tau').
  \end{align*}  
  When $\tau$ is $3k^2$ for some odd $k$, the exact same argument still goes through, except the rational function $p$ is replaced by $q$, which is some other (still fixed and explicit) rational function.  Both $p$ and $q$ are written out on pages 118 and 119 of \cite{Masser1975}, but we will write them here for completeness\footnote{The reader may note a strange-looking asymmetry in $p$ and $q$, namely the $\beta_{0,1}$ in the denominator.  Why not $\beta_{1,0}$?  Masser in fact proves in his work that $\beta_{0,1}=\beta_{1,0}$, so really there is no asymmetry.}.
  \[p(j,\beta_{i,k})=\dfrac{9j(\beta_{2,0}-\beta_{1,1}+\beta_{0,2})}{\beta_{0,1}}+\dfrac{3(7j-6912)}{2(j-1728)}.\]
  \[q(j,\beta_{i,k})=\dfrac{9j(\beta_{4,0}-\beta_{3,1}+\beta_{2,2}-\beta_{1,3}+\beta_{0,4})}{\beta_{0,1}}+\dfrac{3(7j-6912)}{2(j-1728)}.\]  
  In either case we get $\sigma\psi(\tau)=\psi(\tau')$.  Since $\chi^*=r(j,\psi)$ for a $\mathbb{Q}$-rational function $r$, we get $\sigma(\chi^*(\tau))=\chi^*(\tau')$ as required.
 \end{proof}
 
\begin{corollary}\label{cor:StrongGaloisOrbits}
 Let $K$ be a number field.  There are positive constants $c, \delta >0$ with the following property.   Let $\tau\in\mathbb{D}$ be a quadratic point of discriminant $D$.  Then there are $\gg D^\delta$ distinct quadratic points $\tau'\in\mathbb{D}$, of height at most $cD$, such that $(j(\tau'),\chi^*(\tau'))$ is a Galois conjugate, over $K$, of the point $(j(\tau),\chi^*(\tau))$.
 \end{corollary}
 \begin{proof}
  For quadratic points $\sigma\in\mathbb{D}$, let $H(\sigma)$ be the height of $\sigma$ and $D(\sigma)$ the discriminant.  
  It is known that the number of distinct Galois conjugates of $j(\sigma)$ over $\mathbb{Q}$ is bounded from below by a positive power of $D(\sigma)$.  This follows from the Siegel lower bound \cite{Siegel1935} for class numbers of quadratic fields.  See Pila \cite{Pila2011} for more details.  
  
  Since $[K:\mathbb{Q}]$ is a fixed constant, the number of Galois conjugates of $j(\tau)$ over $K$ is therefore $\gg D^\delta$.  Each Galois conjugate $\theta_i$ of $j(\tau)$ over $K$ yields a distinct $\tau_i\in\mathbb{D}$, such that $\theta_i(j(\tau))=j(\tau_i)$.  Moreover, $D(\tau_i)=D(\tau)=D$.  
  
  By work of Pila \cite{Pila2011}, there is a constant $c$ such that, for any $\sigma\in\mathbb{D}$,
  \[H(\sigma)\leq cD(\sigma).\]
  Hence each $\tau_i$ has $H(\tau_i)\leq cD$.
    
  Finally, by \ref{propn:GaloisOrbits}, we have $(j(\tau_i),\chi^*(\tau_i))=(\theta_i(j(\tau)),\theta_i(\chi^*(\tau)))$. 
 \end{proof}
Corollary \ref{cor:StrongGaloisOrbits} gives us exactly the lower bound we need to work with the Pila-Wilkie theorem.  Shortly we will use this bound and the Pila-Wilkie theorem \ref{thrm:PilaWilkie} to prove our main theorem, \ref{thrm:AOforChiStar2}.  First we have a proposition demonstrating the ideas in the simplest case; it also serves as the base case for an inductive argument we use in \ref{thrm:AOforChiStar2}.

\begin{proposition}[Andr\'e-Oort for $(j,\chi^*)$, in 2 Dimensions]\label{propn:AOforChiStar1}
 Let $C\suq\mathbb{C}^2$ be an irreducible algebraic curve.  Then $C$ contains only finitely many $(j,\chi^*)$-special points.
 \end{proposition}
 \begin{proof}
  Suppose that $C$ contained infinitely many special points.  Since special points are algebraic, this tells us that $C$ can in fact be defined over $\overline{\mathbb{Q}}$ and thus over a number field $K$.
  
  Define a set $Z\suq\mathbb{D}$ by
  \[Z=\{\tau\in\mathbb{D}:(j(\tau),\chi^*(\tau))\in C\}.\]
  Then $Z$ is definable. If it contains an arc of a real algebraic curve, then by \ref{thrm:AxLChiStarAndJ} it must be all of $\uh$, which is impossible since $j$ and $\chi^*$ are algebraically independent.  Hence $Z^\text{alg}$ is empty.  We will show that $Z\setminus Z^\text{alg}=Z$ contains `many' (ie. a positive power of $T$) quadratic points of a given height $T$, contradicting the Pila-Wilkie theorem.
  
  Since $C$ contains infinitely many special points, we have infinitely many distinct quadratic points $\tau\in\mathbb{D}$ with $(j(\tau),\chi^*(\tau))\in C$.  In particular, we can find such a $\tau$ with arbitrarily large discriminant $D$.  Hence by \ref{cor:StrongGaloisOrbits}, there are $\gg D^\delta$ quadratic points $\tau'\in\mathbb{D}$, of height at most $cD$, such that $(j(\tau'),\chi^*(\tau'))$ is a Galois conjugate of $(j(\tau),\chi^*(\tau))$ over $K$.  
  
  Since it is a Galois conjugate of $(j(\tau),\chi^*(\tau))$, we know that \[(j(\tau'),\chi^*(\tau'))\in C,\] hence all of the $\tau'$ lie in $Z$.  So there are $\gg D^\delta$ quadratic points (of height at most $cD$) in $Z=Z\setminus Z^\text{alg}$, which contradicts the Pila-Wilkie Theorem for any $\epsilon < \delta$.
 \end{proof}

In more dimensions, the fundamental ideas for dealing with special points by counting Galois conjugates are exactly the same; we have the following.  (Compare with Theorem 11.2 of \cite{Pila2011}.)

\begin{proposition}\label{propn:InductiveCount}
 Suppose $V\subseteq\mathbb{C}^{2n}$ is a variety defined over a number field $K$.  Write $V^{\text{sp}}$ for the union of all positive-dimensional $(j,\chi^*)$-special subvarieties of $V$.  Suppose that $V^{\text{sp}}$ is a variety.  Then $V\setminus V^{\text{sp}}$ contains only finitely many $(j,\chi^*)$-special points.
 \end{proposition}
 \begin{proof}
  Let $\mathcal{Z}=\pi^{-1}(V)$ and $Z=\mathcal{Z}\cap\mathbb{D}^n$.  Then $Z$ is definable.
  
  The set $\mathcal{Z}^\text{alg}$ consists of $\mathcal{Z}^\text{sp}=\pi^{-1}(V^{\text{sp}})$ as well as possibly some weakly $\uh$-special varieties; but the weakly $\uh$-special varieties can contain no quadratic points.  Hence, if we denote by $N(X,T)$ the number of quadratic points in $X$ up to height $T$, we have
  \[N(Z\setminus Z^\text{sp},T)\leq N(Z\setminus Z^\text{alg},T)\ll_\epsilon T^\epsilon\]
  for any $\epsilon>0$; the last bound coming from the Pila-Wilkie Counting Theorem.  Here $Z^\text{sp}=\mathcal{Z}^\text{sp}\cap\mathbb{D}^n$.
  
  Suppose for a contradiction that $V\setminus V^\text{sp}$ contains infinitely many $(j,\chi^*)$-special points.  Then we can find quadratic points 
  \[u=(\tau_1,\dots,\tau_n)\in Z\setminus Z^\text{sp}\]
  of arbitrarily large discriminant $D$.  By \ref{cor:StrongGaloisOrbits}, there are $\gg D^\delta$ quadratic points $u'\in\mathbb{D}^n$, with height at most $cD$, such that $\pi(u')$ is a Galois conjugate of $\pi(u)$ over $K$.  This gives us $\gg D^\delta$ quadratic points (of height at most $cD$) in $Z\setminus Z^\text{sp}$.  Choosing any $\epsilon <\delta$, we get a contradiction to the Pila-Wilkie theorem for sufficiently large $D$.
 \end{proof}

So we have some control over the special points that can arise in a given variety.  The next step is to deal with the positive-dimensional special subvarieties.
\begin{definition}
 A $\uh$-special (or $(j,\chi^*)$-special, or $j$-special, etc.) variety is called \emph{basic} if it has no constant factors.  That is, if the set $S_0$, from the definition of a special variety, is empty.
\end{definition}
Every weakly $\uh$-special variety $S$ arises as the product of a basic $\uh$-special variety $B$ with some number of constant factors $q_i$ (if all the $q_i$ are quadratic points then $S$ is special).  When this happens, we say that $S$ is the \emph{translate} of $B$ by the factors $q_i$.  The following lemma tells us which basic special varieties have translates lying in the preimage of a given variety $V$.

\begin{lemma}\label{lma:BasicPreimages}
 Let $V\subseteq\mathbb{C}^{2n}$ be a variety and define $\mathcal{Z}=\pi^{-1}(V)$.  There is a finite collection $\mathcal{B}$ of basic $\uh$-special varieties with the property that every maximal, positive-dimensional, weakly $\uh$-special subvariety of $\mathcal{Z}$ is a translate of $\gamma B$, for some $B\in\mathcal{B}$ and $\gamma\in\slz^k$.
\end{lemma}
\begin{proof}
  This is identical to Proposition 10.2 of \cite{Pila2011}.  In the presence of \ref{thrm:AxLChiStarAndJ}, the proof carries over exactly.
 \end{proof}

Finally, we combine \ref{propn:InductiveCount} and \ref{lma:BasicPreimages} in an inductive argument to prove our main theorem.

\begin{theorem}[Andr\'e-Oort for $(j,\chi^*)$]\label{thrm:AOforChiStar2}
 Let $V\suq\mathbb{C}^{2n}$ be a variety.  Then $V$ contains only finitely many maximal $(j,\chi^*)$-special subvarieties.
 \end{theorem}
 \begin{proof}
  There is a subvariety $\tilde V\suq V$, defined over $\overline{\mathbb{Q}}$, containing all the algebraic points of $V$.  So we may assume that $V$ is defined over $\overline{\mathbb{Q}}$ (and thus over a number field $K$).
  
  We will proceed by induction on $n$.  The base case is \ref{propn:AOforChiStar1}.  The conclusion holds by \ref{propn:InductiveCount} if $V^\text{sp}$ is variety.  So it is sufficient to prove that $V^\text{sp}$ is a variety, under the assumption that \ref{thrm:AOforChiStar2} holds for $m<n$.
 
  By \ref{lma:BasicPreimages}, there are finitely many basic $\uh$-special varieties, $B\in \mathcal{B}$, such that every maximal $\uh$-special subvariety of $\pi^{-1}(V)$ is a translate of some $\gamma B$.  A maximal $(j,\chi^*)$-special subvariety of $V$ is the Zariski closure of $\pi(S)$, for some maximal $\uh$-special subvariety $S\suq\pi^{-1}(V)$.  Therefore any maximal $(j,\chi^*)$-special subvariety of $V$ is the translate (by some special points $(j(\tau_i),\chi^*(\tau_i))$) of one of a finite collection $\mathcal{C}$ of basic $(j,\chi^*)$-special varieties.  (The twists by elements of $\slz$ have no effect since $j$ and $\chi^*$ are modular.)
  
  So it is enough to show that, given some basic special $C\in\mathcal{C}$, there are only finitely many translates of $C$ which are maximal $(j,\chi^*)$-special subvarieties of $V$.  Such a $C$ will be a subvariety of $\mathbb{C}^{2k}$ for some $k$.
  
  The possible translates of $C$ are elements of $\mathbb{C}^{2(n-k)}$, namely the set of points\footnote{We are being slightly lax with our labelling of coordinates here.  The constant factors by which we translate our basic varieties can be in any of the pairs of coordinates in $(\mathbb{C}^{2})^n$.  Since there are only finitely many ways to reorder the coordinates, no issues will arise from allowing the translations to take place in any of the coordinates.}
  \begin{multline*}V'=\bigl\{(j_1,\chi_1,\dots,j_{n-k},\chi_{n-k}):\text{the translate of }C\text{ by }\\(j_1,\chi_1,\dots,j_{n-k},\chi_{n-k})\text{ is contained in }V\bigr\}.\end{multline*}
  This is an algebraic subvariety of $\mathbb{C}^{2(n-k)}$.  The translates of $C$ which yield special subvarieties of $V$ are the $(j,\chi^*)$-special points of $V'$.  The translates which yield \emph{maximal} special subvarieties are the $(j,\chi^*)$-special points of $V'\setminus(V')^\text{sp}$.  By our inductive assumption, there are only finitely many such points.  Thus $V^\text{sp}$, which consists of finitely many translates of the finitely many basic special varieties in $\mathcal{C}$, is a variety.  So we can conclude by \ref{propn:InductiveCount}.
 \end{proof}

\begin{corollary}[Andr\'e-Oort for $\chi^*$]
 Let $V\suq\mathbb{C}^n$ be a variety.  Then $V$ contains only finitely many maximal $\chi^*$-special subvarieties.
\end{corollary}
\begin{proof}
  Consider a variety $V'\suq \mathbb{C}^{2n}$, defined as
  \[V'=\{(J_1,X_1,\dots,J_n,X_n):(X_1,\dots,X_n)\in V\}.\]
  Given a maximal $\chi^*$-special subvariety $S$ of $V$, there is a corresponding $(j,\chi^*)$-special subvariety $S'\suq V'$, such that the projection of $S'$ onto the $X_i$ coordinates (which correspond to $\chi^*$) is $S$.  By \ref{thrm:AOforChiStar2}, it is enough to show that $S'$ is a maximal $(j,\chi^*)$-special subvariety of $V'$.
  
  Indeed, if $S'$ were contained in a $(j,\chi^*)$-special subvariety $T\suq V'$, with $\dim T >\dim S'$, then by the definition of $(j,\chi^*)$-special varieties, there must be a condition on a $\chi^*$-coordinate which is relaxed in going from $S'$ to $T$.  Hence the projection of $T$ onto the $X_i$ coordinates would be a $\chi^*$-special subvariety of $V$ strictly containing $S$.  Contradiction.
 \end{proof}
\bigskip
\bibliographystyle{../../../bib/scabbrv}
\bibliography{../../../bib/thebib}

\end{document}